\documentclass[doublespacing]{elsart}

 \usepackage[latin1]{inputenc}

\usepackage{natbib,amssymb,enumerate,amsmath,bbm,ifthen}

\newcommand{\wrt}{with respect to}

\newcommand{\iid}{i.i.d.}
\newcommand{\as}{a.s.}

\newcommand{\ie}{i.e.}
\newcommand{\eg}{e.g.}

\newcommand{\eqsp}{\;}
\newcommand{\lleb}{\lambda^{\mathrm{Leb}}}
\newcommand{\nset}{\mathbb{N}}
\newcommand{\rset}{\mathbb{R}}
\newcommand{\mcf}{\mathcal{F}}

\newcommand{\tvnorm}[1]{\ensuremath{\left\|#1\right\|_{\mathrm{TV}}}}

\def\tf{f_\star}
\def\tW{W_\star}
\def\tV{V_\star}
\def\tU{U_\star}

\newcommand{\xprod}{\bar{x}}
\newcommand{\Xprod}{\bar{X}}
\newcommand{\fprod}{\bar{f}}

\newcommand{\Pprod}{\bar{P}}
\newcommand{\Qprod}{\bar{Q}}
\newcommand{\Vprod}{\bar{V}}
\newcommand{\gprod}{\bar{g}}
\newcommand{\lambdaprod}{\bar{\lambda}}

\newcommand{\eqdef}{\ensuremath{\stackrel{\mathrm{def}}{=}}}

\def\rme{\mathrm{e}}

\def\PP{\mathbb{P}}
\def\tphi{\phi_\star}
\def\tsigma{\sigma_\star}
\def\tbeta{\beta_\star}
\def\tPP{\PP_\star}
\def\tPE{\PE_\star}

\def\PE{\mathbb{E}}

\def\1{\mathbbm{1}}

\def\Aset{\mathsf{A}}

\def\Cset{\mathsf{C}}
\def\tCset{\Cset_\star}

\newcommand{\Csetprod}{{\bar{\Cset}}}
\def\Dset{\mathsf{D}}

\def\Kset{\mathsf{K}}
\def\Xset{\mathsf{X}}
\def\Zset{\mathsf{Z}}
\def\Xsetprod{\bar{\Xset}}

\def\Xsigma{\mathcal{X}}
\def\Zsigma{\mathcal{Z}}
\def\Xsigmaprod{\bar{\mathcal{X}}}
\def\xb{\mathbf{x}}

\def\JointKernel{T}
\def\tJointKernel{T_\star}

\def\tb{b_\star}

\def\th{h_\star}
\def\tsigma{\sigma_\star}
\def\tlambda{\lambda_\star}

\newcommand{\Xinit}{\ensuremath{\nu}}
\newcommand{\tXinit}{\ensuremath{\nu_\star}}

\newcommand{\Xproc}{\ensuremath{\{X_k\}_{k\geq 0}}}

\newcommand{\Yset}{\ensuremath{\mathsf{Y}}}
\newcommand{\Ysigma}{\ensuremath{\mathcal{Y}}}
\newcommand{\Yproc}{\ensuremath{\{Y_k\}_{k\geq 0}}}

\newcommand{\XYproc}{\ensuremath{\{X_k,Y_k\}_{k\geq 0}}}

\newcommand{\chunk}[4][]%
{\ifthenelse{\equal{#1}{}}{\ensuremath{{#2}_{#3:#4}}}{\ensuremath{#2^#1}_{#3:#4}}
}

\newcommand{\Q}{\ensuremath{Q}}
\newcommand{\G}{\ensuremath{G}}
\newcommand{\tG}{\ensuremath{G}_\star}
\newcommand{\tQ}{\ensuremath{Q_\star}}
\newcommand{\StatDistrib}{\pi}
\newcommand{\tStatDistrib}{\pi_\star}

\newcommand{\filt}[2][]
{
\ifthenelse{\equal{#1}{}}{\ensuremath{\phi_{#2}}}{\ensuremath{\phi_{#1,#2}}}
}
\newcommand{\pred}[3][]
{
\ifthenelse{\equal{#1}{}}{\ensuremath{\phi_{#2|#3}}}{\ensuremath{\phi_{#1,#2|#3}}}%
}
\newcommand{\post}[3][]
{
\ifthenelse{\equal{#1}{}}{\ensuremath{\phi_{#2|#3}}}{\ensuremath{\phi_{#1,#2|#3}}}%
}
\newcommand{\logl}[2][]%
{%
\ifthenelse{\equal{#1}{}}{\ensuremath{\ell_{#2}}}{\ensuremath{\ell_{#1,#2}}}%
}
\newcommand{\lhood}[2][]%
{%
\ifthenelse{\equal{#1}{}}{\ensuremath{\mathrm{L}_{#2}}}{\ensuremath{\mathrm{L}_{#1,#2}}}%
}

\newcounter{hyp}

\def\ie{\textit{i.e.}}

\newcommand{\Afunc}{\ensuremath{a}}
\newcommand{\Bfunc}{\ensuremath{b}}
\newcommand{\U}{\ensuremath{\zeta}}
\newcommand{\Uproc}{\ensuremath{\{\U_k\}_{k\geq 0}}}
\newcommand{\V}{\ensuremath{\varepsilon}}
\newcommand{\Vproc}{\ensuremath{\{\V_k\}_{k\geq 0}}}

\journal{Stochastic Process. Appl.}

\begin{document}

\begin{frontmatter}
  
  \title{Forgetting of the initial distribution for Hidden Markov Models\thanksref{anr}}

\author[cmapx]{R. Douc},
\ead{douc@cmapx.polytechnique.fr}
\author[ltci]{G. Fort~\corauthref{cor}},
\corauth[cor]{Corresponding author.}
\ead{gfort@tsi.enst.fr}
\author[ltci]{E. Moulines}, 
\ead{moulines@tsi.enst.fr}
\author[paris6]{P. Priouret} 
\ead{priouret@ccr.jussieu.fr}

\address[cmapx]{CMAP, Ecole Polytechnique, 91128 Palaiseau, France.}
\address[ltci]{LTCI, CNRS-GET/Télécom Paris, 46 rue Barrault, 75634 Paris Cedex
  13, France.}  \address[paris6]{LPMA, Univ. P.M. Curie, Boîte courrier
  188, 75252 Paris Cedex 05, France.}  \thanks[anr]{This work was partly
  supported by the National Research Agency (ANR) under the program
  ``ANR-05-BLAN-0299''.}

\begin{abstract}
  The forgetting of the initial distribution for discrete Hidden Markov Models
  (HMM) is addressed: a new set of conditions is proposed, to establish the
  forgetting property of the filter, at a polynomial and geometric rate.  Both
  a pathwise-type convergence of the total variation distance of the filter
  started from two different initial distributions, and a convergence in
  expectation are considered. The results are illustrated using different HMM
  of interest: the dynamic tobit model, the non-linear state space model and
  the stochastic volatility model.
\end{abstract}

\end{frontmatter}

\section{Definition and notations}
A Hidden Markov Model (HMM) is a doubly stochastic process with an underlying
Markov chain that is not directly observable. More specifically, let $\Xset$
and $\Yset$ be two spaces equipped with a countably generated $\sigma$-fields
$\Xsigma$ and $\Ysigma$; denote by $\Q$ and $\G$ respectively, a Markov
transition kernel on $(\Xset,\Xsigma)$ and a transition kernel from
$(\Xset,\Xsigma)$ to $(\Yset,\Ysigma)$.  Consider the Markov transition kernel
defined for any $(x,y) \in \Xset \times \Yset$ and $C \in \Xsigma \otimes
\Ysigma$ by
\begin{equation}
\label{eq:JointChainHMM}
\JointKernel \left[(x,y), C\right] \eqdef \Q \otimes \G [(x,y),C] = \iint \Q(x,dx') \, \G(x',dy') \1_C(x',y') \eqsp.
\end{equation}
We consider $\XYproc$ the Markov chain with transition kernel $\JointKernel$
and initial distribution $\Xinit \otimes \G(C) \eqdef \iint \Xinit(dx) \G(x,dy)
\1_C(x,y)$, where $\Xinit$ is a probability measure on $(\Xset,\Xsigma)$.  We
assume that the chain $\Xproc$ is not observable (hence the name
\emph{hidden}).  The model is said to be partially dominated if there exists a
measure $\mu$ on $(\Yset,\Ysigma)$ such that for all $x \in \Xset$,
$\G(x,\cdot)$ is absolutely continuous \wrt\ $\mu$: in such case, the joint
transition kernel $\JointKernel$ can be written as
\begin{equation}
\label{eq:JointChain:part_dominatedHMM}
\JointKernel\left[(x,y), C\right] = \iint \Q(x,dx') g(x',y')\, \1_C(x',y')\mu(dy') \eqsp, \quad  C \in \Xsigma \otimes \Ysigma \eqsp,
\end{equation}
where $g(x,\cdot)= \frac{d G(x,\cdot)}{d \mu}$ denotes the Radon-Nikodym
derivative of $G(x,\cdot)$ \wrt\ $\mu$. To follow the usage in the filtering
literature, $g(x,\cdot)$ is referred to as the \emph{likelihood} of the
observation.  An example of such type of dependence is $X_{k+1}=
\Afunc(X_k,\U_{k+1})$ and $Y_k= \Bfunc(X_k,\V_k)$, where $\Uproc$ and $\Vproc$
are i.i.d. sequences of random variables, and $\Uproc$, $\Vproc$ and $X_0$ are
independent. The most elementary example is the so-called linear Gaussian state
space model (LGSSM) where $a$ and $b$ are linear and $\{\U_k,\V_k\}_{k \geq 0}$
are \iid\ standard Gaussian. We denote by $\filt[\Xinit]{n}[\chunk{y}{0}{n}]$
the distribution of the hidden state $X_n$ conditionally on the observations
$\chunk{y}{0}{n} \eqdef [y_0, \dots, y_n]$, which is given by
\begin{multline}
  \label{eq:filtering-distribution-1}
  \filt[\Xinit]{n}[\chunk{y}{0}{n}](A) \eqdef \frac{\Xinit \left[ g(\cdot,y_0) \Q g(\cdot,y_1) \Q \dots \Q g(\cdot,y_n) \1_A \right]}{\Xinit \left[ g(\cdot,y_0) \Q g(\cdot,y_1) \Q \dots \Q g(\cdot,y_n) \right]} \\
  = \frac{\int_{\Xset^{n+1}} \Xinit(dx_0) g(x_0,y_1) \prod_{i=1}^n
    \Q(x_{i-1},dx_i) g(x_i,y_i) \1_A(x_n)}{\int_{\Xset^{n+1}} \Xinit(dx_0)
    g(x_0,y_1) \prod_{i=1}^n \Q(x_{i-1},dx_i) g(x_i,y_i)} \eqsp,
\end{multline}
where $\Q f(x)= \Q(x,f) \eqdef \int \Q(x,dx') f(x')$, for any function $f \in
\mathbb{B}_+(\Xset)$ the set of non-negative functions $f : \Xset \to \rset$,
such that $f$ is $\Xsigma/\mathcal{B}(\rset)$ measurable, with
$\mathcal{B}(\rset)$ the Borel $\sigma$-algebra.

In practice the model is rarely known exactly and so suboptimal filters are
constructed by replacing the unknown transition kernel, likelihood function and
initial distribution by suitable approximations.

The choice of these quantities plays a key role both when studying the
convergence of sequential Monte Carlo methods or when analysing the asymptotic
behaviour of the maximum likelihood estimator (see \eg\ \cite{delmoral:2004} or
\cite{cappe:moulines:ryden:2005} and the references therein).

The simplest problem assumes that the transitions are known, so that the only error in the
filter is due to a wrong initial condition.  A typical question is to ask
whether $\filt[\Xinit]{n}[\chunk{y}{0}{n}]$ and
$\filt[\Xinit']{n}[\chunk{y}{0}{n}]$ are close (in some sense) for large values
of $n$, and two different choices of the initial distribution $\Xinit$ and
$\Xinit'$.

The forgetting property of the initial condition of the optimal filter in
nonlinear state space models has attracted many research efforts and it would
be a formidable task to give credit to every contributors. The purpose of the
short presentation of the existing results below is mainly to allow comparison
of assumptions and results presented in this contributions \wrt\ those
previously reported in the literature. The first result in this direction has
been obtained by \cite{ocone:pardoux:1996}, who established $L_p$-type
convergence of the optimal filter initialised with the wrong initial condition
to the filter initialised with the true initial distribution (assuming that the
transition kernels are known); however, their proof does not provide a rate of
convergence.  A new approach based on the Hilbert projective metric has later
been introduced in \cite{atar:zeitouni:1997} to obtain the exponential
stability of the optimal filter \wrt\ its initial condition. However their
results were based on stringent \emph{mixing} conditions for the transition
kernels; these conditions state that there exist positive constants
$\epsilon_-$ and $\epsilon_+$ and a probability measure $\lambda$ on
$(\Xset,\Xsigma)$ such that for $f \in \mathbb{B}_+(\Xset)$,
\begin{equation}
\label{eq:mixing-condition}
\epsilon_- \lambda(f) \leq \Q(x,f) \leq \epsilon_+ \lambda(f) \eqsp, \quad \text{for any $x \in \Xset$} \eqsp.
\end{equation}
This condition in particular implies that the chain is uniformly geometrically
ergodic.  Similar results were obtained independently by
\cite{delmoral:guionnet:1998} using the Dobrushin ergodicity coefficient (see
\cite{delmoral:ledoux:miclo:2003} for further refinements under this
assumption). The mixing condition has later been weakened by
\cite{chigansky:lipster:2004}, under the assumption that the kernel $Q$ is
positive recurrent and is dominated by some reference measure $\lambda$:
\[
\sup_{(x,x') \in \Xset \times \Xset} q(x,x') < \infty \quad \text{and} \quad \int \mathrm{ess inf} q(x,x') \StatDistrib(x) \lambda(dx) > 0 \eqsp,
\]
where $q(x,\cdot)= \frac{d \Q(x,\cdot)}{d\lambda}$, $\mathrm{ess inf}$ is the essential infimum \wrt\ $\lambda$ and
$\StatDistrib d \lambda$ is the stationary distribution of the chain $\Q$ . If the upper
bound is reasonable, the lower bound is restrictive in many applications and
fails to be satisfied \eg\ for the linear state space Gaussian model.

In \cite{legland:oudjane:2003}, the stability of the optimal filter is studied
for a class of kernels referred to as \emph{pseudo-mixing}. The definition of
pseudo-mixing kernel is adapted to the case where the state space is $\Xset=
\rset^d$, equipped with the Borel sigma-field $\Xsigma$.  A kernel $\Q$ on
$(\Xset,\Xsigma)$ is \emph{pseudo-mixing} if for any compact set $\Cset$ with a
diameter $d$ large enough, there exist positive constants $\epsilon_-(d) >0$
and $\epsilon_+(d) > 0$ and a  measure $\lambda_\Cset$ (which may
be chosen to be finite without loss of generality) such that
\begin{equation}
\label{eq:pseudo-mixing-kernel}
\epsilon_-(d) \lambda_\Cset(A)\leq \Q(x,A)\leq \epsilon_+(d) \lambda_\Cset(A) \eqsp, \quad \text{for any $x \in \Cset$, $A \in \Xsigma$}
\end{equation}
This condition implies that for any $(x',x'') \in \Cset \times \Cset$,
$$
\frac{\epsilon_-(d)}{\epsilon_+(d)} < \mathrm{essinf}_{x \in \Xset} q(x',x)/q(x'',x)\leq  \mathrm{esssup}_{x \in \Xset} q(x',x)/q(x'',x) \leq \frac{\epsilon_+(d)}{\epsilon_-(d)} \eqsp,
$$
where $q(x,\cdot) \eqdef d \Q(x,\cdot)/ d \lambda_\Cset$, and $\mathrm{esssup}$ and $\mathrm{essinf}$
denote the essential supremum and infimum \wrt\ $\lambda_\Cset$.  This
condition is obviously more general than \eqref{eq:mixing-condition}, but still
it is not satisfied in the linear Gaussian case (see \cite[Example
4.3]{legland:oudjane:2003}).

Several attempts have been made to establish the stability conditions under the
so-called \emph{small} noise condition. The first result in this direction has
been obtained by \cite{atar:zeitouni:1997} (in continuous time) who considered
an ergodic diffusion process with constant diffusion coefficient and linear
observations: when the variance of the observation noise is sufficiently small,
\cite{atar:zeitouni:1997} established that the filter is exponentially stable.
Small noise conditions also appeared (in a discrete time setting) in
\cite{budhiraja:ocone:1999} and \cite{oudjane:rubenthaler:2005}. These results
do not allow to consider the linear Gaussian state space model with arbitrary
noise variance.

A very significant step has been achieved by
\cite{kleptsyna:veretennikov:2007}, who considered the filtering problem of
Markov chain $\Xproc$ with values in $\Xset= \rset^d$ filtered from
observations $\Yproc$ in $\Yset= \rset^\ell$,
\begin{equation}
\label{eq:NLGSSM}
\begin{cases}
  X_{k+1} = X_k + b(X_k) + \sigma(X_k) \U_k \eqsp, \\
  Y_k = h(X_k) + \beta \V_k \eqsp.
\end{cases}
\end{equation}
Here $\{ (\U_k,\V_k) \}_{k \geq 0} $ is a \iid\ sequence of standard Gaussian
random vectors in $\rset^{d+\ell}$, $b(\cdot)$ is a $d$-dimensional vector
function, $\sigma(\cdot)$ a $d \times d$-matrix function, $h(\cdot)$ is a
$\ell$-dimensional vector-function and $\beta > 0$. The author established,
under appropriate conditions on $b$, $h$ and $\sigma$, that the optimal filter
forgets the initial conditions; these conditions cover (with some restrictions)
the linear Gaussian state space model.

In this contribution, we will propose a new set of conditions to establish the
forgetting property of the filter, which are more general than those proposed
in \cite{kleptsyna:veretennikov:2007}.  In theorem
\ref{thm:almost-sure-convergence}, a pathwise-type convergence of the total
variation distance of the filter started from two different initial
distributions is established, which is shown to hold almost surely w.r.t.  the
probability distribution of the observation process $\{Y_k\}_k$. Then, in
Theorem~\ref{thm:convergence-moment}, the convergence of the expectation of
this total variation distance is shown, under more stringent conditions.  The
results are shown to hold under rather weak conditions on the observation
process $\{Y_k\}_k$ which do not necessarily entail that the observations are
from an HMM.

The paper is organised as followed. In section
\ref{sec:AssumptionsMainResults}, we introduce the assumptions and state the
main results. In section~\ref{sec:ApplicationsHMM}, we give sufficient
conditions for Theorems~\ref{thm:almost-sure-convergence} and
\ref{thm:convergence-moment} to hold, when $\{Y_k\}_k$ is an HMM process,
assuming that the transition kernel and the likelihood function might be
different from those used in the definition of the filter. In section
\ref{sec:Examples}, we illustrate the use of our assumptions on several
examples with unbounded state spaces. The proofs are given in sections
\ref{sec:ProofofTheoremsAlmostSureConvergence} and
\ref{sec:ProofofpropositionsBoundConditionInitiale}.

\section{Assumptions and Main results}
\label{sec:AssumptionsMainResults}
We say that a set $\Cset \in \Xsigma$ satisfies the \emph{local Doeblin}
property (for short, $\Cset$ is a LD-set), if there exists a measure
$\lambda_\Cset$ and constants $\epsilon^-_\Cset > 0$ and $\epsilon^+_\Cset > 0$
such that, $\lambda_\Cset(\Cset) > 0$ and for any $A \in \Xsigma$,
\begin{equation}
\label{eq:definition-LD-set}
\epsilon^{-}_\Cset \lambda_\Cset (A \cap \Cset)\leq \Q(x, A \cap \Cset )\leq \epsilon^{+}_\Cset \lambda_\Cset(A \cap \Cset) \eqsp, \quad
\text{for all $x \in \Cset$} \eqsp.
\end{equation}
Locally Doeblin sets share some similarities with $1$-small set in the theory of Markov chains over general
state spaces (see \cite[chapter 5]{meyn:tweedie:1993}).  Recall that a set
$\Cset$ is $1$-small if there exists a measure
$\tilde{\lambda}_\Cset$ and $\tilde{\epsilon}_\Cset > 0$, such that $\tilde{\lambda}_{\Cset}(\Cset) > 0$, and for all $x
\in \Cset$ and $A \in \Xsigma$, $\Q(x,A \cap \Cset) \geq \tilde{\epsilon}_\Cset \tilde{\lambda}_\Cset(A \cap \Cset)$. In particular, a locally Doeblin set is $1$-small
with $\tilde{\epsilon}_\Cset= \epsilon^-_\Cset$ and $\tilde{\lambda}_\Cset= \lambda_\Cset$.
The main difference stems from the fact that we impose both a lower and an
upper bound, and we impose that the minorizing and the majorizing measure are
the same.

Compared to the pseudo-mixing condition \eqref{eq:pseudo-mixing-kernel}, the
local Doeblin property involves the trace of the Markov kernel $Q$ on $\Cset$
and thus happens to be much less restrictive. In particular, on the contrary to
the pseudo-mixing condition, it can be easily checked that for the kernel
associated to the linear Gaussian state space model, every bounded Borel set
$\Cset$ is locally Doeblin.

Let $V$ be a positive function $V: \Xset \to [1,\infty)$ and $A \in \Xsigma$ be
a set.  Define:
\begin{equation}
\label{eq:definition-Upsilon}
\Upsilon_\Aset(y) \eqdef \sup_{x \in \Aset} g(x,y) QV(x)/V(x) \eqsp.
\end{equation}
Consider the following assumptions:
\begin{enumerate}[(H1)]
\item \label{assum:likelihood-not-zero} For any $(x,y) \in \Xset \times \Yset$, $g(x,y) > 0$.
\item \label{assum:likelihood-drift} There exist a set $\Kset \subseteq \Yset$
  and a function $V: \Xset \to [1,\infty)$ such that for any $\eta > 0$, one may choose a LD-set $\Cset \in \Xsigma$
  satisfying
\[
\Upsilon_{\Cset^c}(y) \leq \eta \; \Upsilon_{\Xset}(y) \eqsp, \quad \text{for all $y \in \Kset$.}
\]
\end{enumerate}
Assumption (H\ref{assum:likelihood-not-zero}) can be relaxed, but this
assumption simplifies the statements of the results and the proofs. The case
where the likelihood may vanish will be considered in a companion paper.
Assumption (H\ref{assum:likelihood-drift}) involves both the likelihood
function and the drift function.  It is satisfied for example if there exists a
set $\Kset$ such that for all $\eta > 0$, one can choose a LD-set $\Cset$ so
that
\begin{equation}
\label{eq:condition-likelihood}
\sup_{x \in \Cset^c} g(x,y) < \eta \sup_{x \in \Xset} g(x,y) \eqsp, \quad \text{for all $y \in \Kset$,}
\end{equation}
in which case the previous assumption is satisfied with $V \equiv 1$. When
$\Xset = \rset^d$, this situation occurs for example when the compact sets are
locally Doeblin and $\lim_{|x| \to \infty} \sup_{y \in \Kset} g(x,y)= 0$. As a
simple illustration, this last property is satisfied for
$Y_k=h(X_k)+\epsilon_k$ with $\lim_{|x|\to \infty} |h(x)| =\infty$ and
$\{\epsilon_k\}_k$ are i.i.d.random variables (independent of $\{X_k\}_k$) with
a density $g$ which satisfies $\lim_{|x|\to \infty} g(x)=0$. More complex
models satisfying (H\ref{assum:likelihood-drift}) are considered in
Section~\ref{sec:Examples}.

When \eqref{eq:condition-likelihood} is not satisfied, assumption
(H\ref{assum:likelihood-drift}) can still be fulfilled if for all $y \in
\Yset$, $\sup_{x \in \Xset} g(x,y) < \infty$, $\sup_{\Xset} \Q V/V < \infty$
for some function $V: \Xset \to [1,\infty)$, and for all $\eta>0$, there exists
a LD-set $\Cset$ such that $\sup_{\Cset^c} \Q V/V \leq \eta$. As a simple
illustration, this situation occurs for example with $X_{k+1}=\phi X_k + \sigma
\zeta_k$ where $|\phi|<1,\ \sigma>0$ and $\{\zeta_k\}_k$ a family of iid
standard Gaussian vectors. More details are provided in Section
\ref{sec:Examples}.

For any LD-set $\Dset$ and $\Xinit$ a probability measure on $(\Xset,\Xsigma)$
define:
\begin{align}
\label{eq:definition-Phi}
&\Phi_{\Xinit,\Dset}(y_0,y_1) \eqdef \Xinit[g(\cdot,y_0) \Q g(\cdot,y_1) \1_{\Dset}] \eqsp, \\
\label{eq:definition-Psi}
&\Psi_\Dset(y) \eqdef \lambda_\Dset ( g(\cdot,y) \1_\Dset ) \eqsp.
\end{align}

We denote by $(\Omega,\mathcal{A})$ a measurable space, and we let $\Yproc$ be
a stochastic process on $(\Omega,\mathcal{A})$ which takes values in
$(\Yset,\Ysigma)$ but \emph{which is not necessarily the observation of an
  HMM}. For any probability measure $\Xinit$ and any $n \in \nset$, the
filtering distribution $\filt[\Xinit]{n}[\chunk{Y}{0}{n}]$ (defined in
\eqref{eq:filtering-distribution-1}) is a measure-valued random variable on
$(\Omega,\mathcal{A})$.
\begin{thm}
\label{thm:almost-sure-convergence}
Assume (H\ref{assum:likelihood-not-zero}-\ref{assum:likelihood-drift}) and let $\tPP$ be a probability measure on $(\Omega,\mathcal{A})$. Assume
in addition that for some LD-set $\Dset$ and some constants $M > 0$ and $\gamma
\in (0,1)$,
\begin{align}
\label{eq:condition-nb-visites-K}
&\liminf_{n \to \infty} n^{-1} \sum_{i=0}^n \1_{\Kset}(Y_i) \geq (1+\gamma)/2 \eqsp,  && \quad \tPP-\as \\
\label{eq:condition-Upsilon}
&\limsup_{n \to \infty} n^{-1} \sum_{i=0}^n \log \Upsilon_{\Xset} (Y_i) < M \eqsp, && \quad \tPP-\as \\
\label{eq:condition-Psi_D}
&\liminf_{n \to \infty} n^{-1} \sum_{i=2}^n \log \Psi_\Dset (Y_i) > - M \eqsp,
&& \quad \tPP-\as
\end{align}
where $\Upsilon_{\Xset}$ and $\Psi_\Dset$ are defined in
\eqref{eq:definition-Upsilon} and \eqref{eq:definition-Psi}, respectively.
Then, for any initial distributions $\Xinit$ and $\Xinit'$ on $(\Xset,
\Xsigma)$ such that $\Xinit(V) + \Xinit'(V) < \infty$, $\Xinit Q \1_\Dset > 0$
and $\Xinit' Q \1_\Dset > 0$, there exists a positive constant $c$ such that,
\begin{equation}
\label{eq:geometric-ergodicity-as}
\limsup_{n \to \infty} n^{-1} \log \tvnorm{\filt[\Xinit]{n}[\chunk{Y}{0}{n}] - \filt[\Xinit']{n}[\chunk{Y}{0}{n}]} < -c  \eqsp, \quad \tPP-\as
\end{equation}
\end{thm}

\begin{rem}
  We stress that it is not necessary to assume that $\Yproc$ is the observation
  of an HMM $\XYproc$. Conditions \eqref{eq:condition-Upsilon} and
  \eqref{eq:condition-Psi_D} can be verified for example under a variety of
  weak dependence conditions, the only requirement being basically to be able
  to prove a LLN (see for example \cite{dedecker:doukhan:2003}). This is of
  interest because in many applications, the HMM model is not correctly
  specified, but it is still of interest to establish the forgetting properties
  of the filtering distribution \wrt\ the initial distribution.
\end{rem}
We will now state a statement allowing to control the expectation of the total variation distance.
\begin{thm}
\label{thm:convergence-moment}
Assume (H\ref{assum:likelihood-drift}).  Let $\Dset$ be a LD-set.  Then, for
any $M_i > 0$, $i=0,1,2$, and $\gamma \in (0,1)$, there exist $\beta \in (0,1)$
such that, for any given initial distributions $\Xinit$ and $\Xinit'$ on
$(\Xset, \Xsigma)$ and all $n$,
\begin{multline}
\label{eq:convergence-moment}
\tPE \left( \tvnorm{\filt[\Xinit]{n}[\chunk{Y}{0}{n}] - \filt[\Xinit']{n}[\chunk{Y}{0}{n}]} \right) \\
\leq \beta^n \; \left[ 1 + \Xinit(V) \Xinit'(V) \right] + r_0(\Xinit,n) + r_0(\Xinit',n) + \sum_{i=1}^3 r_i(n) \eqsp
\end{multline}
where the sequences $\{ r_0(\Xinit,n) \}_{n \geq 0}$ and $\{r_i(n)\}_{n \geq
  0}$, $i=1,2,3$ are defined by
\begin{align}
\label{eq:definition-r0}
&r_0(\Xinit,n) \eqdef \tPP \left( \log \Phi_{\Xinit,\Dset}(Y_0,Y_1) \leq -M_0 n \right) \eqsp, \\
\label{eq:definition-r1}
&r_1(n) \eqdef \tPP \left( \sum_{i=0}^n \log \Upsilon_{\Xset}(Y_i) \geq M_1 n \right) \eqsp, \\
\label{eq:definition-r2}
&r_2(n) \eqdef \tPP \left( \sum_{i=0}^n  \log \Psi_\Dset(Y_i)  \leq  - M_2 n\right) \eqsp, \\
\label{eq:definition-r3}
&r_3(n) \eqdef \tPP \left( n^{-1} \sum_{i=1}^n \1_{\Kset}(Y_i) \leq (1+\gamma)/2 \right)   \eqsp.
\end{align}
\end{thm}
\section{Applications to HMM}
\label{sec:ApplicationsHMM}
We will now discuss conditions upon which \eqref{eq:condition-Upsilon} and
\eqref{eq:condition-Psi_D} hold (Propositions~\ref{prop:check-assumption-LLN}
to \ref{prop:drift-condition-LLN}) and upon which the right hand side in
(\ref{eq:convergence-moment}) vanishes
(Proposition~\ref{prop:bound-condition-initiale} to
Corollary~\ref{cor:drift-condition-MLP}). To that goal, we assume that $\Yproc$
is the observation of an HMM $\XYproc$ with Markov kernel $\tJointKernel= \tQ
\otimes \tG$, where $\tQ$ is a transition kernel on $(\Xset,\Xsigma)$ and $\tG$
is a Markov kernel from $(\Xset,\Xsigma)$ to $(\Yset,\Ysigma)$, and initial
distribution $\tXinit$ on $(\Xset,\Xsigma)$.

Recall that a kernel $P$ on a general state space $(\Zset,\Zsigma)$ is
phi-irreducible and (strongly) aperiodic if there exists a $\sigma$-finite
measure $\varphi$ on $(\Zset,\Zsigma)$, such that, for any $A \in \Zsigma$
satisfying $\varphi(A) > 0$ and any initial condition $x$, $P^n(x;A) > 0$, for
all $n$ sufficiently large. A set $\Cset \in \Zsigma$ is called \emph{petite}
for the Markov kernel $P$ if for some probability measure $m$ on $\nset$, with
finite mean sampling time (which can always be done without loss of
generality~\cite[Proposition 5.5.6]{meyn:tweedie:1993})
\[
\sum_{n=0}^\infty m(n) P^n(x,A) \geq \epsilon_\Cset^- \lambda_\Cset(A) \eqsp, \quad \text{for all $x \in \Cset$, $A \in \Zsigma$},
\]
where $\lambda_\Cset$ is a measure on $(\Zset,\Zsigma)$ satisfying $\lambda_\Cset(\Cset) > 0$ and $\epsilon_\Cset^- >0$. We
denote by $\PP_\Xinit^P$ and $\PE_\Xinit^P$ the probability distribution and the expectation on the canonical probability space
$(\Zset^\nset,\Zsigma^{\otimes \nset})$ associated to the Markov chain with transition kernel $P$ and initial distribution $\Xinit$.

We first state sufficient conditions for $\tJointKernel$ to be an aperiodic
positive Harris chain (see definitions and main properties in \cite[Chapters 10
\& 13]{meyn:tweedie:1993} and \cite[Chapter 14]{cappe:moulines:ryden:2005}) and
for the law of large numbers to hold for the Markov chain with kernel
$\tJointKernel$.

\begin{prop}
\label{prop:check-assumption-LLN}
Assume that $\tQ$ is an aperiodic, positive Harris Markov kernel with
stationary distribution $\tStatDistrib$.  Then, the kernel $\tJointKernel$
defined by 
\[
\tJointKernel[(x,y),A] \eqdef \iint \tQ(x,dx') \tG(x',dy') \1_A(x',y') \eqsp,
\quad A \in \Xsigma \otimes \Ysigma \eqsp,
\]
is an aperiodic positive Harris Markov kernel with stationary distribution
$\tStatDistrib \otimes \tG$. In addition, for any initial distribution
$\tXinit$ on $(\Xset,\Xsigma)$, and any function $\varphi \in
\mathbb{B}_+(\Xset \times \Yset)$ satisfying $\tStatDistrib \otimes \tG
(\varphi) < \infty$,
\begin{equation}
\label{eq:check-assumption}
n^{-1} \sum_{i=0}^n \varphi(X_i,Y_i) \to \tStatDistrib \otimes \tG (\varphi) \quad \PP^{\tJointKernel}_{\tXinit \otimes \tG}-\as\
\end{equation}
\end{prop}
\begin{cor}
\label{cor:Finiteness-StationaryDistribution}
If $\tStatDistrib \otimes \tG \left( \log \Upsilon_\Xset \right)_+ < \infty$
(resp. $\tStatDistrib \otimes \tG \left( \log \Psi_\Dset \right)_- < \infty$),
then, condition \eqref{eq:condition-Upsilon} (resp. \eqref{eq:condition-Psi_D})
is satisfied with $\tPP \eqdef \PP^{\tJointKernel}_{\tXinit \otimes \tG}$.
\end{cor}
In many problems of interest, it is not straightforward to establish that the
chain is positive Harris; in addition, the distribution $\tStatDistrib$ is not
known explicitly making the conditions of Corollary
\ref{cor:Finiteness-StationaryDistribution} difficult to check.  It is often
interesting to apply the following result which is a direct consequence of the
$f$-norm ergodic theorem and the law of large numbers for positive Harris chain
(see for example \cite[Theorems 14.0.1, 17.0.1]{meyn:tweedie:1993}).
\begin{prop}
\label{prop:drift-condition-LLN}
Let $\tf \geq 1$ be a function on $\Xset$. Assume that $\tQ$ is a
phi-irreducible Markov kernel and that there exist a petite set $\tCset$, a
function $\tV: \Xset \to [1,\infty)$, and a constant $\tb$ satisfying
\begin{equation}
\label{eq:drift-condition-additive}
\tQ \tV(x) \leq \tV(x) - \tf(x) + \tb \1_{\tCset}(x)\eqsp.
\end{equation}
Then, the kernel $\tQ$ is positive Harris with invariant probability
$\tStatDistrib$ and $\tStatDistrib(\tf) < +\infty$.  Let $\varphi \in
\mathbb{B}_+(\Xset \times \Yset)$ be a function such that
\begin{equation}
\label{eq:f-dominate}
\sup_{x \in \Xset} \tf^{-1}(x)  \tG \left(x, \varphi(x,\cdot) \right) < \infty,
\end{equation}
Then, $\tStatDistrib \otimes \tG (\varphi) < \infty$.
\end{prop}

We now derive conditions to compute a bound for $\{r_0(\Xinit,n)\}_{n \geq 0}$.
\begin{prop}
\label{prop:bound-condition-initiale}
Assume (H\ref{assum:likelihood-not-zero}-\ref{assum:likelihood-drift}) and that
the drift function $V$ defined in (H\ref{assum:likelihood-drift}) satisfies
$\sup_\Xset V^{-1} \Q V < \infty$.
\begin{enumerate}[(i)]
\item If for some $p \geq 1$,
\begin{equation}
\label{eq:condition-initiale-Lp}
\sup_{i=0,1} \sup_\Xset V^{-1} \tPE [\log g(\cdot,Y_i) ]_-^p < \infty \eqsp,
\end{equation}
then, there exists a constant $C$ such that, for any initial probability
measure $\Xinit$ on $(\Xset,\Xsigma)$ such that $\nu Q \1_D>0$ and all $n \geq
0$, $r_0(\Xinit,n) \leq C n^{-p} \Xinit(V)$.
\item If for some positive $\lambda$,
\begin{equation}
\label{eq:condition-initiale-exponentiale}
\sup_{i=0,1} \sup_\Xset V^{-1} \tPE \left(\exp( \lambda [\log g(\cdot,Y_i) ]_-) \right)< \infty \eqsp,
\end{equation}
then there exist positive constants $C,\delta > 0$, such that for any initial
probability measure $\Xinit$ on $(\Xset,\Xsigma)$ such that $\nu Q \1_D>0$, and
all $n \geq 0$, $r_0(\Xinit,n) \leq C \rme^{-\delta n} \Xinit(V)$.
\end{enumerate}
\end{prop}

To determine the rate of convergence of the sequences $\{r_i(n) \}_{n \geq 0}$
to zero, $i=1,2,3$, it is required to use deviation inequalities for partial
sums of the observations $\Yproc$. There are a variety of techniques to prove
such results, depending on the type of assumptions which are available. If
polynomial rates are enough, then one can apply the standard Markov inequality
together with the Marcinkiewicz-Siegmund inequality; see for example
\cite{dedecker:doukhan:2003} or \cite{fort:moulines:2003}.
\begin{prop}
\label{prop:drift-condition-Polynomial}
Assume that
\begin{enumerate}[(i)]
\item $\tQ$ is  aperiodic and positive Harris Markov kernel with stationary distribution $\tStatDistrib$.
\item \label{item:driftcondition} There exist a petite set $\tCset$ and functions
  $\tU,\tV,\tW: \Xset \to [1,\infty)$ and a constant $\tb$ satisfying $\tStatDistrib(\tW) < \infty$ and
  \begin{align*}
  & \tQ \tV \leq \tV - \tU + \tb \1_{\tCset} \eqsp, \\
  & \tQ \tW \leq \tW - \tV + \tb \1_{\tCset}
  \end{align*}
\end{enumerate}
Let $p \geq 1$. There exists a constant $C < \infty$ such that for any function
$\varphi$ on $(\Yset, \Ysigma)$ satisfying $\sup_\Xset \tU^{-1}
\tG(\cdot,|\varphi|^p) < \infty$ and $\sup_\Xset \tU^{-1} \tV^{1-1/p} \;
\tG(\cdot,|\varphi|) < \infty$, and for any initial probability distribution
$\tXinit$ on $(\Xset, \Xsigma)$, and any $\delta > 0$,
\begin{align*}
  & \PP^{\tJointKernel}_{\tXinit \otimes \tG} \left[ \sum_{i=1}^n \left\{
      \varphi(Y_i) - \tStatDistrib \otimes \tG(\varphi) \right\} \geq \delta n
  \right] \leq C \delta^{-p} n^{-(p/2 \vee 1)} \tXinit(\tW)\eqsp,
\end{align*}
 \end{prop}
 \begin{cor}
   If there exists $p \geq 1$ such that
\[
\sup_\Xset \tf^{-1} \tG(\cdot,|\log \Upsilon_\Xset|^p) < \infty \eqsp, \qquad
\sup_{\Xset} \tf^{-1} \tV^{1-1/p} \tG(\cdot,|\log \Upsilon_\Xset|) < \infty
\eqsp,
\]
and
\[
\sup_\Xset \tf^{-1} \tG(\cdot,|\log \Psi_\Dset|^p) < \infty \eqsp, \qquad \sup_{\Xset}
\tf^{-1} \tV^{1-1/p} \tG(\cdot,|\log \Psi_\Dset|) < \infty \eqsp,
\]
then there exist finite constants $C, M_i$, $i=1,2,3$ such that
\[
r_i(n) \leq C n^{-(p/2 \vee 1)} \tXinit(\tW) \eqsp.
\]
 \end{cor}

 If we wish to establish that the sequences $\{r_i(n)\}_{n \geq 0}$ decreases
 to zero exponentially fast, we might for example use the multiplicative
 ergodic theorem \cite[Theorem 1.2]{kontoyiannis:meyn:2005} to bound an
 exponential moment of the partial sum, and then use the Markov inequality.
 This will require to check the multiplicative analog of the additive drift
 condition \eqref{eq:drift-condition-additive}.

Some additional definitions are needed.  Let $W: \Xset \to (0,\infty)$
be a function.  We say that the function $W$ is \emph{unbounded} if $\sup_\Xset W = + \infty$. We define by
$\mathcal{G}_{W}$ the set of functions whose growth at infinity is lower than $W$, \ie\ $F$ belongs to $\mathcal{G}_{W}$ if and only if
\begin{equation}
\label{eq:definition-mathcalG}
\sup_\Xset \left( |F| - W \right) < \infty \eqsp.
\end{equation}
\begin{prop}
\label{prop:drift-condition-MLP}
Let $\tW$ be an unbounded function $\tW: \Xset \to (0,\infty)$ and that the
level sets $\{ \tW \leq r \}$ are petite. Assume that $\tQ$ is phi-irreducible
and that there exist a function $\tV: \Xset \to [1,\infty)$, and constant $\tb
< \infty$ such that
\begin{equation}
\log \left(  \tV^{-1} \tQ \tV \right) \leq -\tW + \tb \eqsp. \label{eq:drift-condition-MLP}
\end{equation}
Then, $\tQ$ is positive Harris with a unique invariant probability distribution $\tStatDistrib$, satisfying
$\tStatDistrib(\tV) < \infty$. Let $\varphi$ be a non-negative function. If for some $\tlambda > 0$,
\begin{equation}
\label{eq:condition-growth-Upsilon-varphi}
\log \left[ \tG \left( \cdot, \rme^{\tlambda \varphi} \right) \right] \in \mathcal{G}_{\tW} \eqsp,
\end{equation}
there exists a constant $M > 0$ such that, for any initial distribution $\tXinit$ satisfying $\tXinit \left( \tV \right) < \infty$,
\begin{equation}
\label{eq:borne-exponentielle:varphi}
\limsup_{n \to \infty} n^{-1} \log \PP_{\tXinit \otimes \tG}^{\tJointKernel} \left( \sum_{i=0}^n \varphi(Y_i) \geq M n \right) < 0 \eqsp.
\end{equation}
\end{prop}
\begin{cor}
\label{cor:drift-condition-MLP}
Assume that for some $\tlambda > 0$, 
\[
\log \left[ \tG \left( \cdot, \rme^{\tlambda [\log \Upsilon_\Xset]_+} \right)
\right] \in \mathcal{G}_{\tW} \qquad \log \left[ \tG \left( \cdot,
    \rme^{\tlambda [\log \Psi_\Dset]_-} \right) \right] \in \mathcal{G}_{\tW} \eqsp.
\]
Then, there exist constants $M_i$, $i=1,2$ such that $\limsup_{n \to \infty}
n^{-1} \log r_i(n) < 0$, where $\{ r_i(n)\}_{n \geq 0}$ are defined in
\eqref{eq:definition-r1} and \eqref{eq:definition-r2}.
\end{cor}

\section{Examples}
\label{sec:Examples}
In this section, we illustrate our results using different models of interest.
\subsection{The dynamic tobit  model }
\label{subsec:dynamicTobit}
The tobit model is simply the time series extension of the standard univariate
tobit model and so the univariate hidden process is only observed when it is
positive (\cite{manrique:shephard:1998} and \cite{andrieu:doucet:2002}):
\begin{equation}
\label{eq:tobit-model}
\begin{cases}
X_{k+1}= \phi X_{k}+\sigma \U_{k}\eqsp, \\
Y_k=\max(X_k+\beta \V_k,0)\eqsp,
\end{cases}
\end{equation}
where $\{(\U_k,\V_k) \}_{k \geq 0} $ is a sequence of \iid\ standard Gaussian
vectors, and $|\phi| < 1$, $\sigma > 0$ and $\beta > 0$.  Here $\Xset= \rset$,
$\Yset= \rset_+$ and $\Xsigma$ and $\Ysigma$ are the corresponding Borel
$\sigma$-algebra.  The model is partially dominated (see
\eqref{eq:JointChain:part_dominatedHMM}) \wrt\ the dominating measure
$\delta_0+ \lleb$, where $\lleb$ is the Lebesgue measure and $\delta_0$ is the
Dirac mass at zero. The transition kernels $Q_{\phi,\sigma}$ and the likelihood
$g_\beta$ are respectively given by:
\begin{align}
\label{eq:transition-kernel-tobit}
&Q_{\phi,\sigma}(x,A)= \left( 2 \pi \sigma^{2} \right)^{-1/2} \int \exp \left[ -(1/2 \sigma^{2}) (x'- \phi x)^2 \right] \1_A(x') \lleb(dx') \eqsp, \\
\label{eq:likelihood-tobit} \nonumber
&g_\beta(x,y)= \1 \{y=0 \} (2 \pi \beta^{2})^{-1/2} \int_{x}^{\infty} \exp\left[ -(1/2 \beta^{2}) v^2 \right] \lleb(dv) \\
&\hspace{100pt} + \1 \{ y > 0 \} (2 \pi \beta^{2})^{-1/2} \exp \left[ -(1/2 \beta^{2}) (y-x)^2 \right] \eqsp.
\end{align}
We denote  $Q= Q_{\phi,\sigma}$ and $g= g_\beta$. \\
We assume that $\Yproc$ are the observations of a tobit model
\eqref{eq:tobit-model} with initial distribution $\tXinit$ and 'parameters'
$\tphi$, $\tsigma$, $\tbeta$ (which may be different from $\phi$, $\sigma$,
$\beta$) satisfying $|\tphi| < 1$, $\tsigma > 0$ and $\tbeta > 0$. We denote by
$\tQ= Q_{\tphi,\tsigma}$, $\tG(x,\cdot)= g_{\tbeta}(x,\cdot) \lleb$ and $\tPE=
\PE_{\tXinit \otimes \tG}^{\tJointKernel}$, where $\tJointKernel= \tQ \otimes
\tG$.

\subsubsection{Assumptions H\ref{assum:likelihood-not-zero} and H\ref{assum:likelihood-drift} }It is easily seen that any bounded Borel set $\Cset \subset \{x, 0 \leq |x|
\leq C \}$ satisfies the local Doeblin property \eqref{eq:definition-LD-set},
with $\lambda_\Cset(\cdot)= (2C)^{-1} \lleb \left(\1_\Cset \cdot \right)$.
Assumption (H\ref{assum:likelihood-not-zero}) is trivially satisfied. To check
(H\ref{assum:likelihood-drift}), we set $\Kset=\Yset$ and $V(x) = e^{c|x|}$ for
some $c>0$.  The function $V^{-1} Q V $ is locally bounded and $\lim_{|x| \to
  \infty} V^{-1}(x) Q V(x) = 0$. Therefore, since $\sup_{\Xset \times \Yset}
g(x,y) \leq 1\vee (2 \pi \beta^2)^{-1/2}$, for any $\eta > 0$ one may choose a
constant $C>0$ large enough so that $ \Upsilon_{\Cset^c}(y) \leq \eta
\Upsilon_{\Xset}(y)$, where $\Cset \eqdef \{0 \leq |x| \leq C\}$ and
$\Upsilon_A$ is defined in \eqref{eq:definition-Upsilon}.  Therefore,
(H\ref{assum:likelihood-drift}) is satisfied.

\subsubsection{Application of Theorem~\ref{thm:almost-sure-convergence}}
We now check conditions (\ref{eq:condition-nb-visites-K}) to
(\ref{eq:condition-Psi_D}) of Theorem~\ref{thm:almost-sure-convergence}.
Conditions (\ref{eq:condition-nb-visites-K}) and \eqref{eq:condition-Upsilon}
are obvious since $\Kset = \Yset$ and $\sup_{\Yset} \Upsilon_{\Xset} < \infty$.
We now check \eqref{eq:condition-Psi_D} with $\Dset= \{ 0 \leq |x| \leq D\}$
and $\lambda_\Dset(\cdot)= (2D)^{-1} \lleb \left( \1_\Dset \cdot \right)$ where
the constant $D$ is an arbitrary positive constant.  $\tQ(x,dy)$ is a Gaussian
density with mean $\tphi x$ and standard deviation $\tsigma$. Using standard
arguments, $\tQ$ is aperiodic, positive Harris with invariant distribution
$\tStatDistrib$ which is a centered gaussian distribution with variance
$\tsigma^2/(1-\tphi^2)$, and any compact set is petite. By the Jensen
inequality, $\log \lambda_\Dset( g(\cdot,y) \1_\Dset)= \log
\lambda_\Dset(g(\cdot,y)) \geq \lambda_\Dset \left( \log g(\cdot,y) \right)$,
which implies
\begin{multline}
\label{eq:lowerbound_PsiD:tobit}
\log \Psi_\Dset(y) = \log \lambda_\Dset(g(\cdot,y))  \geq \1 \{y=0\} \log \left\{ (2 \pi \beta^2)^{-1/2} \int_{D}^{\infty} \rme^{-v^2/2\beta^2} \lleb(dv) \right\} \\
+ \1 \{ y > 0 \} \left\{ -(1/2) \log (2 \pi \beta^2) - (12 D \beta^2)^{-1}
  \left( (D+y)^3 + (D-y)^3\right) \right\} \eqsp,
\end{multline}
so that $\tStatDistrib \otimes \tG (\left[ \log \Psi_\Dset \right]_- ) <
\infty$.  Corollary \ref{cor:Finiteness-StationaryDistribution} implies
\eqref{eq:condition-Psi_D}.  Combining the results above, Theorem
\ref{thm:almost-sure-convergence} therefore applies showing that
\eqref{eq:geometric-ergodicity-as} holds for any probability $\Xinit$ and
$\Xinit'$ such that $\int \Xinit(dx) \rme^{c |x|} + \int \Xinit'(dx) \rme^{c
  |x|} < \infty$ for some $c > 0$.

\subsubsection{Application of Theorem~\ref{thm:convergence-moment}}
We now consider the convergence of the expectation of the total variation
distance at a polynomial rate.  For all $p \geq 1$, there exists a constant $C$
such that, for any $i \in \{0,1 \}$, $\tPE[ Y_i^{2p}] \leq C ( 1 +
\tPE[X_i^{2p}] )$ which is finite since $\{X_i\}$ is Gaussian.  Therefore,
\begin{equation}
\label{eq:condition:TOBIT-1}
\sup_\Xset (1 + |x|^2)^{-p} \tPE \left[ \log g(x,Y_i) \right]_-^p < \infty \eqsp,
\end{equation}
which implies (\ref{eq:condition-initiale-Lp}) since $V(x) = \exp(c|x|)$. By
Proposition \ref{prop:bound-condition-initiale}, there exists a constant $C$
such that for any probability measure $\Xinit$ such that $\Xinit(V) < \infty$,
$r_0(\Xinit,n) \leq C n^{-p} \Xinit(V)$.  Since $\sup_\Yset \Upsilon_\Xset <
\infty$, we may choose $M_1 > 0$ such that $M_1 > \sup_{\Yset} \log
\Upsilon_\Xset$; for this choice, $r_1(n) \equiv 0$, where $\{r_1(n)\}_{n \geq
  0}$ is defined in \eqref{eq:definition-r1}. Since $\Kset = \Yset$, $r_3(n)
\equiv 0$, where $\{ r_3(n)\}_{n \geq 0}$ is defined in
\eqref{eq:definition-r3}.  We now consider $\{ r_2(n)\}_{n \geq 0}$ and apply
Proposition \ref{prop:drift-condition-Polynomial}. To that goal, we further
assume that there exists $ p_\star \geq 1$ such that $\tXinit(|x|^{3p_\star+1})
< \infty$. It is easily seen that the drift condition
(\ref{eq:drift-condition-additive}) is satisfied with $\tV(x)= 1 +
|x|^{3p_\star}$ and $\tf \sim |x|^{3p_\star-1}$; furthermore, upon noting that
$[\log \Psi_\Dset(y)]_- \sim |y|^2$, we have
\[
\sup_\Xset \tf^{-1} \tG(\cdot, [\log \Psi_\Dset(y)]_-^{p_\star}) < \infty, \qquad
\sup_\Xset \tf^{-1} \tV^{1-1/p_\star} \tG(\cdot, [\log \Psi_\Dset(y)]_-) < \infty
\eqsp,
\]
thus proving $\limsup_{n \to \infty} n^{-(p_\star/2 \vee 1)} r_2(n)= 0$.
Therefore, by Theorem~\ref{thm:convergence-moment}, the expectation $\tPE
\left( \tvnorm{\filt[\Xinit]{n}[\chunk{Y}{0}{n}] -
    \filt[\Xinit']{n}[\chunk{Y}{0}{n}]} \right)$ goes to zero at the rate
$n^{p_\star/2 \vee 1}$ for any initial distributions $\Xinit, \Xinit'$ such
that
$\int \{\Xinit(dx) +  \Xinit'(dx) \} \exp(c|x|)  < +\infty$. \\
The exponential decay can be proved similarly under the assumption that for
some $c>0$, $\int \tXinit(dx) \exp( c |x|) < +\infty$; details are omitted.

\subsection{Non-linear State-Space models}
\label{subsec:NLGSSM}
We consider the model \eqref{eq:NLGSSM} borrowed from
\cite{kleptsyna:veretennikov:2007}. Assume that $\beta>0$,
\begin{list}{}{}
\item[NLG($b,h$)] \label{NLG:locally-bounded}
The functions $b$ and $h$ are locally bounded and
\begin{equation}
\label{eq:NLG:go-to-infty}
\lim_{|x| \to \infty} \left( |x + b(x) | - |x| \right) = - \infty  \eqsp.
\end{equation}
\item[NLG($\sigma$)] \label{NLG:variance} The noise variance is non-degenerated,
\begin{equation}
\label{eq:NLG:variance}
0 <  \inf_{x \in \rset^d} \inf_{\{\lambda \in \rset^d, |\lambda|= 1\}} \lambda^T\sigma(x)\sigma^T(x) \lambda
\leq  \sup_{x \in \rset^d} \sup_{\{\lambda \in \rset^d,|\lambda|= 1\}} \lambda^T\sigma(x)\sigma^T(x) \lambda < \infty \eqsp.
\end{equation}
\end{list}
The model is partially dominated \wrt\ the Lebesgue measure. The transition
kernel $Q_{b,\sigma}$ and the likelihood $g_{h,\beta}$ are respectively given
by
\begin{align}
  & Q_{b,\sigma}(x,A) = (2 \pi)^{-d/2} |\sigma(x)|^{-1} \int \exp \left( -(1/2)  |x'-x-b(x)|_{\sigma(x)}^2 \right) \1_A(x')  \lleb(dx')  \eqsp, \\
  & g_{h,\beta}(x,y)= (2 \pi \beta^2)^{-\ell/2} \exp( - |y-h(x)|^2 / 2 \beta^2)
  \eqsp,
\end{align}
where $|u|_{\sigma(x)}^2= u^T [\sigma(x) \sigma^T(x)]^{-1} u$. As above, we set
$\Q = Q_{b,\sigma}$ and $g = g_{h,\beta}$.

Assume that $\Yproc$ are the observations of a non-linear Gaussian state space
\eqref{eq:NLGSSM} with initial distribution $\tXinit$ and 'parameters' $\tb$,
$\th$, $\tsigma$ and $\tbeta$.  We assume that $\tbeta > 0$ and that the
functions $\tb$, $\th$ and $\tsigma$ satisfy NLG($\tb,\th$)-NLG($\tsigma$),
respectively, and
\begin{equation}
\label{eq:growth-th}
\limsup_{|x| \to \infty} |x|^{-1} \log |\th(x)| < \infty \eqsp.
\end{equation}
We denote by $\tQ= Q_{\tb,\tsigma}$, $\tG= g_{\th,\tbeta} \lleb$ and $\tPE=
\PE_{\tXinit \otimes \tG}^{\tJointKernel}$ where $\tJointKernel= \tQ \otimes
\tG$.
\subsubsection{Assumptions H\ref{assum:likelihood-not-zero} and H\ref{assum:likelihood-drift}}
Under NLG($b,h$)-NLG($\sigma$), every bounded Borel set in $\rset^d$ is locally
Doeblin in the sense given by \eqref{eq:definition-LD-set}.
(H\ref{assum:likelihood-not-zero}) is trivial. Set $V(x)= \exp(c |x|)$, where
$c$ is a positive constant. The likelihood $g$ is bounded by $(2 \pi
\beta^2)^{-\ell/2}$ and under NLG($\sigma$), there exists a constant $M <
\infty$ such that $V^{-1}(x) Q V(x) \leq M \exp\left[ c ( |x + b(x)| - |x|
  )\right]$.  Therefore, under NLG($b,h$)-NLG($\sigma$), for any $\eta > 0$, we
may choose a constant $C$ large enough such that $\Upsilon_{\Cset^c}(y) \leq
\eta \Upsilon_\Xset(y)$ for any $y \in \Yset$ where $\Cset = \{ x \in \rset^d,
|x| \leq C \}$. Hence, assumption (H\ref{assum:likelihood-drift}) is satisfied
with $\Kset=\Yset$.

\subsubsection{Application of Theorem~\ref{thm:almost-sure-convergence}}
Condition~(\ref{eq:condition-nb-visites-K}) is trivial since $\Kset = \Yset$.
Condition \eqref{eq:condition-Upsilon} is obvious too since $\Upsilon_{\Xset}$
is everywhere bounded. For~(\ref{eq:condition-Psi_D}), let us apply
Corollary~\ref{cor:Finiteness-StationaryDistribution} and
Proposition~\ref{prop:drift-condition-LLN}. $\tQ$ is aperiodic, phi-irreducible
and compact sets are petite. Set $\Dset = \{ x \in \rset^d, |x| \leq D \}$,
where $D > 0$ and define $\lambda_{\Dset}(\cdot) = \lleb(\1_{\Dset} \cdot)/
\lleb(\Dset)$.  Noting that $|y-h(x)|^2 \leq 2 ( |y|^2 + |h(x)|^2)$,
\begin{equation}
\label{eq:bound:log-likelihood}
\left[ \log g(x,y) \right]_- \leq \beta^{-2} |y|^2  + \beta^{-2}  |h(x)|^2 + (\ell /2) \left[ \log(2 \pi \beta^2) \right]_+ \eqsp.
\end{equation}
Since the function $h$ is locally bounded, $\sup_\Dset |h|^2 < \infty$ and
\eqref{eq:bound:log-likelihood} implies that
\begin{equation}
\label{eq:lowerbound-Sum}
\left[ \log \Psi_{\Dset}(y) \right]_- \leq  \lambda_{\Dset}(\left[ \log g(\cdot,y) \right]_-)\leq  \beta^{-2} |y|^2  + \beta^{-2}  \sup_\Dset |h|^2 + (\ell /2) \left[ \log(2 \pi \beta^2) \right]_+ \eqsp.
\end{equation}
We set $V_\star(x) = \rme^{c_\star |x|}$ we may find a compact (and thus
petite) set $\Cset_\star$ and constants $\lambda_\star \in (0,1)$ and $s_\star$
such that $\tQ V_\star \leq \lambda_\star V_\star + s_\star \1_{\Cset_\star}$,
so that (\ref{eq:drift-condition-additive}) is satisfied with $\tf=
(1-\tlambda) \tV$. Hence $\tQ$ is positive Harris-recurrent and
$\tStatDistrib(\tV)<+\infty$. Furthermore, Eq.~\eqref{eq:lowerbound-Sum}
implies that there exists a constant $C < \infty$ such that
\begin{equation}
\label{eq:condition-growth-Glog}
\tG\left(x,\left[ \log \Psi_{\Dset} \right]_-\right) \leq C \left( 1 + |\th(x)|^2  \right) \leq C \left( 1+ \tV(x) \sup_\Xset \tV^{-1} |\th|^2 \right) \eqsp.
\end{equation}
The RHS is finite, provided $c_\star \geq 2 \limsup_{|x| \to \infty} |x|^{-1}
\log |\th(x)|$ which we assume hereafter. Therefore, by
Corollary~\ref{cor:Finiteness-StationaryDistribution} and
Proposition~\ref{prop:drift-condition-LLN}, \ref{thm:almost-sure-convergence}
applies: \eqref{eq:geometric-ergodicity-as} holds for any initial probability
measure such that $\int \rme^{c |x|} \Xinit(dx) + \int \rme^{c |x|} \Xinit'(dx)
< + \infty$ for some $c > 0$.

\subsubsection{Application of Theorem~\ref{thm:convergence-moment}}
We are willing to establish geometric rate of convergence and for that purpose
we will use Proposition \ref{prop:bound-condition-initiale} and Proposition
\ref{prop:drift-condition-MLP}. We set $W(x)=  c\{ |x| - |x+b(x)|\} \vee 1$ and
$\tW(x)= c_\star  \{ |x| - |x + \tb(x)| \} \vee 1$ and assume that
\begin{equation}
\label{eq:condition-growth-th}
|h|^2 \in \mathcal{G}_{W} \quad \text{and} \quad |\th|^2  \in \mathcal{G}_{\tW} \eqsp.
\end{equation}

$\tW$ is unbounded and the level sets are petite for $\tQ$. Furthermore,
$\tV(x) = \rme^{c_\star |x|}$ where $c_\star > 0$ satisfies the multiplicative
drift condition \eqref{eq:drift-condition-MLP}.  Let $\lambda < \beta^2(2
\wedge \tbeta^{-2}) /4$.  Since $\lambda \beta^{-2} < \tbeta^{-2}/4$,
Eq.~\eqref{eq:bound:log-likelihood} implies that there exists a constant $C <
\infty$ such that for any integer $i$,
\[
\tPE \left[ \rme^{\lambda [ \log g (x,Y_i) ]_-} \right] \leq C \tPE \left[
  \rme^{2 \lambda \beta^{-2} |\th(X_i)|^2} \right] \rme^{ \lambda \beta^{-2}
  |h(x)|^2} \eqsp.
\]
Since $\lambda \leq \beta^2/2$, Lemma \ref{lem:exponential-bound} shows that
$\sup_i \tPE \left[ \rme^{2 \lambda \beta^{-2} |\th(X_i)|^2} \right] < \infty$
provided $\tXinit(\tV) < +\infty$ which is henceforth assumed. Therefore,
Proposition \ref{prop:bound-condition-initiale} applies, showing that there
exists $\delta >0$ such that for any probability measure $\Xinit$ such that
$\Xinit(V) < \infty$, $r_0(\Xinit,n) \leq C \rme^{-\delta n} \Xinit(V)$.  As in
Section~\ref{subsec:dynamicTobit}, because $\Upsilon_\Xset$ is bounded, we may
choose $M_1$ large enough so that $r_1(n) \equiv 0$ (see
\eqref{eq:definition-r1}); similarly, since $\Kset = \Yset$, $r_3(n) \equiv 0$.
Eq.~(\ref{eq:lowerbound-Sum}) implies that, for any $\tlambda$ small enough,
$\log \tG\left(\cdot,\rme^{\tlambda \left[ \log \Psi_\Dset \right]_-} \right)
\in \mathcal{G}_{\tW}$.  Proposition \ref{prop:drift-condition-MLP} shows that
$\limsup_{n \to \infty} n^{-1} \log r_2(n) < 0$. Hence
Theorem~\ref{thm:convergence-moment} applies: for any initial distribution
$\Xinit,\Xinit'$ such that $\int \{ \Xinit(dx) + \Xinit'(dx) \}\exp(c |x|) <
+\infty$, $\tPE \left( \tvnorm{\filt[\Xinit]{n}[\chunk{Y}{0}{n}] -
    \filt[\Xinit']{n}[\chunk{Y}{0}{n}]} \right)$ goes to zero at a geometric
rate.

\subsection{Stochastic Volatility Model}
\label{subsection:stovol}
As a final example, we consider the stochastic volatility (SV) model.  In the
canonical model in SV for discrete-time data
\cite{hull:white:1987,jacquier:polson:rossi:1994}, the observations $\Yproc$
are the compounded returns and $ \Xproc$ is the log-volatility, which is
assumed to follow a stationary auto-regression of order $1$, \ie\
\begin{equation}
\label{eq:stochastic-volatility}
\begin{cases}
  X_{k+1} = \phi X_k + \sigma \U_k \eqsp, \\
  Y_k = \beta \exp(X_k/2) \V_k \eqsp,
\end{cases}
\end{equation}
where $\{(\U_k,\V_k)\}_{k \geq 0}$ is a \iid\ sequence of standard Gaussian
vectors, $|\phi|<1$, $\sigma>0$ and $\beta >0$.  Here $\Xset=\Yset=\rset$ and
$\Xsigma$ and $\Ysigma$ are the Borel sigma-fields.  The model is partially
dominated \wrt\ the Lebesgue measure. The transition kernel $Q_{\phi,\sigma}$
and the likelihood $g_\beta$ are respectively given by
\begin{align}
& Q_{\phi,\sigma}(x,A)= (2 \pi \sigma^2)^{-1/2} \int \exp( -1/(2 \sigma^2) (x'-\phi x)^2 \1_A(x') \lleb(dx') \eqsp, \\
& g_\beta(x,y)= (2 \pi \beta^2)^{-1/2} \exp\left( - y^2 \exp(-x)/2\beta^2 - x/2 \right) \eqsp.
\end{align}
We denote $Q= Q_{\phi,\sigma}$ and $g= g_\beta$. \\
We assume that $\Yproc$ are the observations of the stochastic volatility model
\eqref{eq:stochastic-volatility} with initial distribution $\tXinit$ and
parameters $|\tphi| < 1$, $\tsigma > 0$, and $\tbeta > 0$. We denote as above
$\tQ= Q_{\tphi,\tsigma}$, $\tG = g_{\tbeta} \lleb$, $\tJointKernel= \tQ \otimes
\tG$ and $\tPE= \PE_{\tXinit \otimes \tG}^{\tJointKernel}$.

\subsubsection{Assumptions H\ref{assum:likelihood-not-zero} and H\ref{assum:likelihood-drift}}
As in example \ref{subsec:dynamicTobit}, every bounded Borel set is locally
Doeblin in the sense of \eqref{eq:definition-LD-set}.  Assumption
(H\ref{assum:likelihood-not-zero}) is satisfied but the likelihood is not
uniformly bounded over $\Xset \times \Yset$; nevertheless it is easily seen
that $\sup_{x \in \Xset} g(x,y) \leq (2 \pi \rme)^{-1/2} |y|^{-1}$. We set
$\Kset= \rset$ and put $V(x)= \rme^{c|x|}$ where $c$ is positive; as in Example
\ref{subsec:dynamicTobit}, $QV(\cdot)/V(\cdot)$ is locally bounded and
$\lim_{|x| \to \infty} Q V(x) / V(x)= 0$, showing that assumption
(H\ref{assum:likelihood-drift}) is fulfilled.

\subsubsection{Application of Theorem~\ref{thm:almost-sure-convergence}}
The Markov kernel $\tQ$ is positive recurrent, geometrically ergodic and its
stationary distribution $\tStatDistrib$ is Gaussian with mean 0 and variance
$\tsigma^2 / (1-\tphi^2)$. Note that there exists a constant $C<\infty$ such
that for all $y\in \Yset$, $\left[ \log \Upsilon_\Xset (y) \right]_+ \leq
C\left| \log|y| \right|$, which implies that $\tG(x,[\log \Upsilon_\Xset ]_+) <
C+|x|/2$ for some constant $C<\infty$. This implies that $\tStatDistrib \otimes
\tG([\log \Upsilon_\Xset ]_+)<\infty$ and
Corollary~\ref{cor:Finiteness-StationaryDistribution} implies
(\ref{eq:condition-Upsilon}). Set $\Dset= \{ x, |x| \leq D \}$ where $D > 0$
and let $\lambda_\Dset(\cdot)= \lleb( \1_\Dset \cdot)/\lleb(\Dset)$.  By the
Jensen inequality,
\[
\log \Psi_\Dset(y) \geq \lambda_\Dset (\log g(\cdot,y))= -(1/2) \log( 2 \pi
\beta^2) - y^2 \; \mathrm{sh}(D)/[2 \beta^2 D] \eqsp,
\]
showing that there exists a constant $C < \infty$ such that $\left[ \log
  \Psi_\Dset(y) \right]_- \leq C(1+y^2)$. Therefore, $\tG(x,[\log
\Psi_\Dset]_-) \leq C(1 + \beta^2 \rme^x)$.  The conditions of Corollary
\ref{cor:Finiteness-StationaryDistribution} are satisfied, showing that
(\ref{eq:condition-Psi_D}) holds.  As a result,
\eqref{eq:geometric-ergodicity-as} holds for any initial distributions $\Xinit$
and $\Xinit'$ such that $\int \Xinit(dx) \exp(c|x|) + \int \Xinit'(dx) \exp(c|x|) < \infty$.

The problem of computing the convergence rates can be addressed as in the other
examples.

\section{Proof of Theorems \ref{thm:almost-sure-convergence} and \ref{thm:convergence-moment}}
\label{sec:ProofofTheoremsAlmostSureConvergence}
Before proving the main results, some additional definitions are needed.  A
function $\bar f$ defined on $\Xsetprod \eqdef \Xset \times \Xset$ is said to
be {\em symmetric} if for all $(x,x') \in \Xsetprod$, $f(x,x')=f(x',x)$. An
unnormalised transition kernel $P$ on $(\Xsetprod, \Xsigmaprod)$, where
$\Xsigmaprod= \Xsigma \otimes \Xsigma$ is said to be {\em symmetric} if for all
$(x,x')$ in $\Xsetprod$ and any positive symmetric function $f$,
$P\left[(x,x'),f \right]= P\left[(x',x),f \right]$.  For $P$ a Markov kernel on
$(\Xset,\Xsigma)$, we denote by $\Pprod$ the transition kernel on
$(\Xsetprod,\Xsigmaprod)$ defined, for any $(x,x') \in \Xsetprod$ and $A$, $A'
\in \Xsigma$, by
\begin{equation}
\label{eq:product-kernel}
\Pprod[ (x,x'), A \times A']= P(x,A) P(x',A') \eqsp.
\end{equation}

For any $A \in \Xsigma$, and $\Xinit$ and $\Xinit'$ two probability
distributions on $(\Xset,\Xsigma)$ the difference
$\filt[\Xinit]{n}[\chunk{y}{0}{n}](A)- \filt[\Xinit']{n}[\chunk{y}{0}{n}](A)$
may be expressed as
\begin{align}
\label{eq:difference}
& \filt[\Xinit]{n}[\chunk{y}{0}{n}](A) - \filt[\Xinit']{n}[\chunk{y}{0}{n}](A) \\ \nonumber
&  = \frac{\PE_\Xinit^{\Q} \left[ \prod_{i=0}^n g(X_i,y_i) \1_A(X_n)\right]}{\PE^{\Q}_\Xinit\left[ \prod_{i=0}^n g(X_i,y_i)\right]} - \frac{\PE^{\Q}_{\Xinit'}\left[ \prod_{i=0}^n g(X_i,y_i) \1_A(X_n)\right]}{\PE^{\Q}_{\Xinit'}\left[ \prod_{i=0}^n g(X_i,y_i)\right]} \\ \nonumber
&  = \frac{\PE_{\Xinit \otimes \Xinit'}^{\Qprod} \left[ \prod_{i=0}^n \gprod(X_i,X'_i,y_i) \1_A(X_n)\right]
- \PE_{\Xinit' \otimes \Xinit}^{\Qprod} \left[ \prod_{i=0}^n \gprod(X_i,X'_i,y_i)  \1_A(X_n)\right]}
{\PE_{\Xinit}^{\Q} \left[ \prod_{i=0}^n g(X_i,y_i) \right] \PE_{\Xinit'}^{\Q} \left[ \prod_{i=0}^n g(X_i,y_i) \right]} \eqsp,
\end{align}
where $\gprod(x,x',y)= g(x,y) g(x',y)$. The idea of writing the difference
using a pair of independent processes has been apparently introduced in
\cite{budhiraja:ocone:1997}; this approach is central in the work of
\cite{kleptsyna:veretennikov:2007}. We consider separately the numerator and
the denominator of Eq.~\eqref{eq:difference}. For the numerator, the path of
the independent processes is decomposed along the successive visits to $\Cset
\times \Cset$ as done in~\cite{kleptsyna:veretennikov:2007}.

\begin{prop}
\label{prop:numerateur}
Let $\Cset$ be a  LD-set and  $\Xinit$ and $\Xinit'$ be two probability distributions on $(\Xset,\Xsigma)$.  For any integer $n$ and functions $g_i \in \mathbb{B}_+(\Xset)$, $i=0,\dots,n$,
such that $\PE^\Q_{\Xinit}\left[ \prod_{i=0}^n g_i(X_i) \right] < \infty$ and $\PE^\Q_{\Xinit'} \left[ \prod_{i=0}^n g_i(X_i) \right] < \infty$, define
\begin{align}
\label{eq:definition-Delta}
&\Delta_n(\Xinit,\Xinit',\{g_i\}_{i=0}^n) \\ \nonumber & \quad = \sup_{A \in
  \Xsigma}\left|\PE_{\Xinit \otimes \Xinit'}^{\Qprod}\left[\prod_{i=0}^n
    \gprod_i(X_i,X'_i) 1_A(X_n)\right] -\PE_{\Xinit' \otimes
    \Xinit}^{\Qprod}\left[\prod_{i=0}^n \gprod_i(X_i,X'_i)
    1_A(X_n)\right]\right| \eqsp,
\end{align}
where $\gprod_i(x,x')= g_i(x)g_i(x')$. Then,
\begin{equation}
\label{eq:bound-Delta-n}
\Delta_n(\Xinit,\Xinit',\{g_i\}_{i=1}^n) \leq \PE^{\Qprod}_{\Xinit \otimes \Xinit'}\left[ \prod_{i=0}^{n} \gprod_i(X_i,X'_i) \rho_\Cset^{N_{\Cset,n}} \right]\eqsp,
\end{equation}
where $\Qprod$ is defined as in \eqref{eq:product-kernel} and
\begin{align}
\label{eq:definition-NCn}
&N_{\Cset,n} \eqdef \sum_{i=0}^{n-1}\1_{\Cset \times \Cset}(X_i,X'_i) \1_{\Cset \times \Cset}(X_{i+1},X'_{i+1}) \eqsp, \\
& \label{eq:definition-rho} \rho_\Cset\eqdef 1- \left(\epsilon^-_\Cset/
  \epsilon^+_\Cset \right)^2\eqsp.
\end{align}
\end{prop}

\begin{pf}
  Put $\xprod=(x,x')$, $\gprod_i(\xprod)=g_i(x)g_i(x')$, $\Csetprod \eqdef
  \Cset\times \Cset$, and $\lambdaprod_\Csetprod \eqdef\lambda_\Cset \otimes
  \lambda_\Cset$.  We stress that the kernels that will be defined along this
  proof may be unnormalized. Since $\Cset$ is a locally Doeblin set, we have for
  any measurable positive function $\fprod$ on $(\Xsetprod,\Xsigmaprod)$,
\begin{equation}
\label{eq:minorization}
(\epsilon^-_\Cset)^2  \lambdaprod_\Csetprod (\1_{\Csetprod} \fprod)\leq \Qprod(\xprod,\1_{\Csetprod} \fprod )\leq (\epsilon^{+}_\Cset)^2 \lambdaprod_\Csetprod(\1_{\Csetprod} \fprod)\eqsp, \quad  \text{for all $\xprod \in \Csetprod$} \eqsp.
\end{equation}
Define the unnormalised kernel $\Qprod_0$ and $\Qprod_1$ by
\begin{align}
\label{eq:definition-Qprod0}
&\Qprod_0(\xprod,\fprod)\eqdef\1_{\Csetprod}(\xprod)(\epsilon^-_\Cset)^2 \lambdaprod_\Csetprod(\1_{\Csetprod} \fprod) \\
\label{eq:definition-Qprod1}
&\Qprod_1(\xprod,\fprod)\eqdef \Qprod(\xprod, \fprod)-\1_{\Csetprod}(\xprod)(\epsilon^-_\Cset)^2 \lambdaprod_\Csetprod(\1_{\Csetprod} \fprod)  = \Qprod(\xprod, \fprod)- \Qprod_0(\xprod,\fprod) \eqsp.
\end{align}
Eq.~\eqref{eq:minorization}  implies that, for all $\xprod \in \Csetprod$,
$
0\leq  \Qprod_1(\xprod, \1_{\Csetprod} \fprod) \leq 
\rho_\Cset \Qprod(\xprod, \1_{\Csetprod} \fprod)$.  It then follows using
straightforward algebra that,
\begin{align}
\label{eq:P1Q}
\Qprod_1(\xprod,\fprod) &= \1_{\Csetprod} (\xprod) \Qprod_1(\xprod,
\1_{\Csetprod} \fprod) + \1_{\Csetprod}(\xprod) \Qprod_1(\xprod,
\1_{\Csetprod^c} \fprod) + \1_{\Csetprod^c}(\xprod) \Qprod_1(\xprod,\fprod) \\
\nonumber &\leq \rho_\Cset \1_{\Csetprod} (\xprod) \Qprod(\xprod,
\1_{\Csetprod} \fprod) + \1_{\Csetprod}(\xprod) \Qprod(\xprod, \1_{\Csetprod^c}
\fprod) + \1_{\Csetprod^c}(\xprod) \Qprod(\xprod,\fprod) \\ \nonumber &\leq
\Qprod(\xprod,\rho_\Cset^{\1_{\Csetprod}(\xprod) \1_{\Csetprod}} \fprod) \eqsp.
\end{align}
We write $\Delta_n(\Xinit,\Xinit',\{g_i\}_{i=0}^n) = \sup_{A \in \Xsigma} |\Delta_n(A)|$
where
\begin{multline}
  \Delta_n(A)\eqdef \Xinit \otimes \Xinit' \left(\gprod_0 \Qprod \gprod_1
    \cdots \gprod_{n-1}\Qprod \gprod_{n} \1_{A \times \Xset} \right) \\ -
  \Xinit' \otimes \Xinit \left(\gprod_0 \Qprod \gprod_1 \cdots \gprod_{n-1}
    \Qprod \gprod_{n} \1_{A \times \Xset} \right) \eqsp.
\end{multline}
Note that $\Delta_n(A)$ may be decomposed as $\Delta_n(A)=\sum_{t_{0:n-1} \in
  \{0,1\}^{n}}\Delta(A,t_{0:n-1})$ where
\begin{multline*}
\Delta_n(A,t_{0:n-1})
\eqdef \Xinit \otimes \Xinit' \left(\gprod_0 \Qprod_{t_0} \gprod_1 \cdots \gprod_{n-1}\Qprod_{t_{n-1}} \gprod_{n} \1_{A\times \Xset} \right) \\ -
\Xinit' \otimes  \Xinit \left( \gprod_0 \Qprod_{t_0} \gprod_1 \cdots \gprod_{n-1}\Qprod_{t_{n-1}} \gprod_{n}\1_{A\times \Xset} \right) \eqsp.
\end{multline*}
Note that, for any $t_{0:n-1} \in \{0,1\}^{n}$ and any sets $A,B \in \Xsigma$,
\begin{multline}
\label{eq:symmetry}
\Xinit' \otimes  \Xinit \left( \gprod_0 \Qprod_{t_0} \gprod_1 \cdots \gprod_{n-1}\Qprod_{t_{n-1}} \gprod_{n}\1_{A \times B} \right) \\
= \Xinit \otimes \Xinit' \left( \gprod_0 \Qprod_{t_0} \gprod_1 \cdots
  \gprod_{n-1}\Qprod_{t_{n-1}} \gprod_{n}\1_{B \times A} \right) \eqsp.
\end{multline}
First assume that there
exists an index $i \geq 0$ such that $t_i=0$ then,
\begin{align*}
&\Xinit \otimes \Xinit' \left(\gprod_0 \Qprod_{t_0} \gprod_1 \cdots \gprod_{n-1}\Qprod_{t_{n-1}} \gprod_{n}1_{A\times \Xset} \right)\\
&\quad =\Xinit \otimes \Xinit' \left(\gprod_0 \Qprod_{t_0} \gprod_1 \cdots \Qprod_{t_{i-1}}\gprod_{i}\1_{\Csetprod}\right)\times (\epsilon^{-}_\Cset)^2 \lambdaprod_\Csetprod \left( \1_{\Csetprod} \gprod_{i+1} \Qprod_{t_{i+1}}  \dots \gprod_{n-1}\Qprod_{t_{n-1}} \gprod_{n}1_{A\times \Xset}\right)\\
&\quad =\Xinit' \otimes \Xinit \left( \gprod_0 \Qprod_{t_0} \gprod_1 \cdots \Qprod_{t_{i-1}}\gprod_{i}\1_{\Csetprod}\right)\times (\epsilon^{-}_\Cset)^2 \lambdaprod_\Csetprod \left(\1_{\Csetprod} \gprod_{i+1} \Qprod_{t_{i+1}} \dots \gprod_{n-1}\Qprod_{t_{n-1}} \gprod_{n}1_{A\times \Xset}\right)
\end{align*}
by (\ref{eq:symmetry}). Thus, $\Delta_n(A,t_{0:n-1})$ is equal to $0$ except if
for all $i$, $t_i=1$, and \eqref{eq:symmetry} finally implies
\[
\Delta_n(A)= \Xinit \otimes \Xinit' \left[\gprod_0 \Qprod_{1} \gprod_1 \cdots
  \gprod_{n-1}\Qprod_{1} \gprod_{n}(\1_{A\times \Xset}-\1_{\Xset \times A})
\right]\eqsp.
\]
Using (\ref{eq:P1Q}), we have
\begin{multline*}
  \Delta_n(\Xinit,\Xinit',\{g_i\}_{i=0}^n) \leq \Xinit \otimes \Xinit'
  \left(\gprod_0 \Qprod_{1} \gprod_1 \cdots \gprod_{n-1}\Qprod_{1} \gprod_{n}
  \right) \\ \leq
  \PE^{\Qprod}_{\Xinit \otimes \Xinit'}\left[ \prod_{i=0}^{n}\gprod_i(\Xprod_i)
    \rho_\Cset^{\sum_{i=0}^{n-1}\1_{\Csetprod}(\Xprod_i)\1_{\Csetprod}(\Xprod_{i+1})}\right]\eqsp,
\end{multline*}
where the last equality is straightforward to establish by induction
on $n$.
The proof is completed.
\end{pf}

\begin{rem}
  If the whole state space $\Xset$ is a locally Doeblin, then one may take $\Cset =
  \Xset$ in the previous expression. Since $N_{\Xset,n}= n$,
  \eqref{eq:difference} and the previous proposition therefore imply the
  uniform ergodicity of the filtering distribution, for any initial
  distribution $\Xinit$ and $\Xinit'$, and any sequence $\chunk{y}{0}{n}$,
  $\tvnorm{\filt[\Xinit]{n}[\chunk{y}{0}{n}] -
    \filt[\Xinit']{n}[\chunk{y}{0}{n}]} \leq \rho_\Xset^n$ where $\rho_\Xset
  \eqdef 1 - (\epsilon_\Xset^-/ \epsilon_\Xset^+ )^2$; see
  \cite{atar:zeitouni:1997} and \cite{delmoral:guionnet:2001}.
 \end{rem}
We consider now the denominator of \eqref{eq:difference}. A lower bound for the
denominator has been computed in \cite[Lemma 2.2]{budhiraja:ocone:1999}. This
is obtained by using a change of measure ideas. We use here a more
straightforward argument.
\begin{prop}
\label{prop:denominateur}
For any LD-set $\Cset \in \Xsigma$, $n\geq 1$ and any functions $g_i \in
\mathbb{B}_+(\Xset)$, $i=0,\dots,n$,
\[
\PE^\Q_{\Xinit} \left[ \prod_{i=0}^n g_i(X_i) \right] \geq (\epsilon^-_\Cset)^{n-1} \Xinit(g_0 \Q  g_1 \1_\Cset) \prod_{i=2}^n \lambda_\Cset (g_i \1_\Cset) \eqsp.
\]
\end{prop}
\begin{pf}
The proof follows immediately from
\[
\PE^\Q_{\Xinit} \left[ \prod_{i=0}^n g_i(X_i) \right] \geq \PE^\Q_{\Xinit} \left[ g_0(X_0) \prod_{i=1}^n g_i(X_i) \1_\Cset(X_i) \right] \eqsp,
\]
and the minorization condition (\ref{eq:definition-LD-set}).
\end{pf}
By combining Propositions \ref{prop:numerateur} and \ref{prop:denominateur}, we can obtain an explicit  bound for the total variation distance
$\tvnorm{\filt[\Xinit]{n}[\chunk{y}{0}{n}]- \filt[\Xinit']{n}[\chunk{y}{0}{n}]}$.

\begin{lem}
\label{lem:main-statement}
Let $\beta \in (0,1)$.  Then, for any LD-sets $\Cset \subseteq \Xset$ and
$\Dset \subseteq \Xset$, any initial probability measures $\Xinit$ and
$\Xinit'$, any function $V : \Xset \to [1,\infty)$,
\begin{multline*}
  \tvnorm{\filt[\Xinit]{n}[\chunk{y}{0}{n}] -
    \filt[\Xinit']{n}[\chunk{y}{0}{n}]} \leq \rho_\Cset^{\beta n} \\ +
  \frac{\prod_{i=0}^n \Upsilon_{\Xset}(y_{i}) \; \max_{\mathcal{I} \subset \{0,
      \dots, n\}, |\mathcal{I}|= a_n} \prod_{i \in \mathcal{I}}
    \Upsilon_{\Cset^c}(y_{i}) \prod_{i \not \in \mathcal{I}}
    \Upsilon_{\Xset}(y_{i}) } {(\epsilon^-_\Dset)^{2(n-1)} \Phi_{\Xinit,\Dset}
    (y_0,y_1) \; \Phi_{\Xinit',\Dset}(y_0,y_1) \prod_{i=2}^n \Psi^2_\Dset(y_i)}
  \Xinit(V) \Xinit'(V) \eqsp,
\end{multline*}
where  $a_{n} \eqdef \lfloor n(1-\beta)/2 \rfloor$, $|\mathcal{I}|$ is the cardinal of the set $\mathcal{I}$  and the functions
$\Phi_{\Xinit,\Dset}$ and $\Psi_{\Dset}$ are defined in \eqref{eq:definition-Phi} and \eqref{eq:definition-Psi}, respectively.
\end{lem}
\begin{pf}
Eq.~\eqref{eq:bound-Delta-n} implies that for any $\beta \in (0,1)$,
\begin{multline*}
  \Delta_n(\Xinit,\tilde{\Xinit},\{g_i\}_{i=0}^n) \leq \PE^{\Qprod}_{\Xinit
    \otimes \Xinit'}\left[ \prod_{i=0}^{n} \gprod(\Xprod_i,y_i)
    \rho_\Cset^{N_{\Cset,n}} \1 \{N_{\Cset,n} \geq \beta n\} \right] \\ +
  \PE^{\Qprod}_{\Xinit \otimes \Xinit'}\left[ \prod_{i=0}^{n}
    \gprod(\Xprod_i,y_i) \rho_\Cset^{N_{\Cset,n}} \1 \{N_{\Cset,n} < \beta n\}
  \right] \eqsp.
\end{multline*}
The first term in the RHS is bounded by $\rho_{\Cset}^{\beta n}
\PE^{\Qprod}_{\Xinit \otimes \Xinit'}\left[ \prod_{i=0}^{n}
  \gprod(\Xprod_i,y_i) \right]$. We now consider the second term. For any set
$\Aset \in \Xsigmaprod$, denote by $M_{\Aset,n}$ the number of visits of $\{
\Xprod_k \}_{k \geq 0}$ to the set $\Aset$ before $n$. By Lemma
\ref{lem:denombrement}, the condition $N_{\Cset,n} < \beta n$ implies that
$M_{\Csetprod,n} < n(1+\beta)/2$ and $M_{\Csetprod^c,n} \geq a_n$.  Note that
for any $\xprod \in \Xsetprod$ and $y \in \Yset$,
\begin{equation}
\label{eq:key-inequality}
\gprod(\xprod,y) \Qprod \Vprod(\xprod) \leq 
[A(y)]^{\1_{\Csetprod^c}(\xprod)}  [B (y)]^{\1_{\Csetprod}(\xprod)} \Vprod(\xprod) \eqsp,
\end{equation}
where we have set $\Vprod(\xprod) \eqdef V(x) V(x')$, $A(y) \eqdef \sup_{\xprod
  \in \Csetprod^c} \gprod(\xprod,y) \Vprod^{-1}(\xprod) \Qprod \Vprod(\xprod)
$, and $B(y) \eqdef \sup_{\xprod \in \Xsetprod} \gprod(\xprod,y)
\Vprod^{-1}(\xprod) \Qprod \Vprod(\xprod) $.  Consider the process
\begin{equation}
\label{eq:definition-Vn}
V_0= \Vprod(\Xprod_0), \quad \text{and} \quad V_n \eqdef \left\{ \prod_{i=0}^{n-1} \frac{\gprod( \Xprod_i,y_i)}{[A(y_i)]^{\1_{\Csetprod^c}(\Xprod_i)}  [B(y_i)]^{\1_{\Csetprod}(\Xprod_i)}}  \right\} \, \Vprod(\Xprod_n) \eqsp, n \geq 1 \eqsp,
\end{equation}
where by convention we have set $0/0= 0$ (to deal with cases where either
$A(y)= 0$ or $B(y)= 0$).  The process $\{V_n\}_{n \geq 0}$ is a
$\mcf$-super-martingale, where $\mcf= \{ \mcf_n \}$ is the natural filtration of
the process $\{ \Xprod_k \}_{k \geq 0}$, $\mcf_n \eqdef \sigma (\Xprod_0,
\dots, \Xprod_n)$. Denote by $\tau_{a_n}$ the $a_n$-th return time to the set
$\Csetprod^c$.  On the event $\{ M_{\Csetprod^c,n} \geq a_n \}$, $\tau_{a_n}
\leq n$, using that $A(y) \leq B(y)$
\begin{multline*}
  \prod_{i=0}^n [A(y_i)]^{\1_{\Csetprod^c}(\Xprod_i)} [B(y_i)]^{\1_{\Csetprod}(\Xprod_i)} \leq \\
  \prod_{i=0}^{\tau_n} [A(y_i)]^{\1_{\Csetprod^c}(\Xprod_i)}
  [B(y_i)]^{\1_{\Csetprod}(\Xprod_i)} \prod_{i=\tau_{a_n}+1}^n B(y_i) \leq
  C(\chunk{y}{0}{n})
\end{multline*}
where $C(\chunk{y}{0}{n}) \eqdef \max_{\mathcal{I} \subset \{0,\dots,n\},
  |\mathcal{I}|= a_n} \prod_{i \in \mathcal{I}} A(y_i) \prod_{i \not \in
  \mathcal{I}} B(y_i)$.  Therefore,
\begin{align*}
&  \PE^{\Qprod}_{\Xinit \otimes \Xinit'}\left[ \prod_{i=0}^n \gprod(\Xprod_i,y_i) \1 \{ N_{\Cset,n} < \beta n \} \right] \leq
\PE^{\Qprod}_{\Xinit \otimes \Xinit'}\left[ \prod_{i=0}^n \gprod(\Xprod_i,y_i) \1 \{ M_{\Csetprod^c,n} \geq a_n \} \right] \\
& \quad \leq C(\chunk{y}{0}{n}) \PE^{\Qprod}_{\Xinit \otimes \Xinit'}\left[ \prod_{i=0}^n \frac{\gprod(\Xprod_i,y_i)}{[A(y_i)]^{\1_{\Csetprod^c}(\Xprod_i)} [B(y_i)]^{\1_{\Csetprod}(\Xprod_i)}} \Vprod(\Xprod_{n+1}) \right] \\
&\quad  =C(\chunk{y}{0}{n})  \PE^{\Qprod}_{\Xinit \otimes \Xinit'} [V_{n+1}]\eqsp.
\end{align*}
The super-martingale inequality therefore implies
\[
\PE^{\Qprod}_{\Xinit \otimes \Xinit'}\left[ \prod_{i=0}^n \gprod(\Xprod_i,y_i) \1 \{ N_{\Cset,n} < \beta n \}  \right] \leq C(\chunk{y}{0}{n})  \Xinit(V) \Xinit'(V) \eqsp,
\]
and the proof follows from \eqref{eq:difference} and Proposition
\ref{prop:denominateur}, using that $A(y) \leq \Upsilon_\Xset(y)
\Upsilon_{\Cset^c}(y)$ and $B(y)= \Upsilon_{\Xset}^2(y)$, where $\Upsilon_A(y)$
is defined in \eqref{eq:definition-Upsilon}.
\end{pf}

\begin{cor}
\label{coro:main-statement}
Assume (H\ref{assum:likelihood-drift}). Let $\Dset$ be a LD-set, and $\gamma$ and $\beta$ be constants satisfying $\gamma \in (0,1)$ and $\beta \in (0,\gamma)$. Then, for any $\eta \in (0,1)$ there exists a LD-set $\Cset$ such that, for
any sequence $\chunk{y}{0}{n} \in \Yset^{n+1}$ satisfying $\sum_{i=0}^n \1_{\Kset}(y_i) \geq (1+\gamma)n/2$, any initial probability measures $\Xinit$ and $\Xinit'$, and
any $n \geq 1$,
\begin{multline*}
  \tvnorm{\filt[\Xinit]{n}[\chunk{y}{0}{n}] -
    \filt[\Xinit']{n}[\chunk{y}{0}{n}]} \leq \rho_\Cset^{\beta n} \\ +
  \frac{\eta^{(\gamma-\beta)n/2} \prod_{i=0}^n \Upsilon^2_\Xset(y_i) }
  {(\epsilon^-_\Dset)^{2(n-1)} \Phi_{\Xinit,\Dset} (y_0,y_1) \;
    \Phi_{\Xinit',\Dset}(y_0,y_1) \prod_{i=2}^n \Psi^2_\Dset(y_i)} \; \Xinit(V)
  \Xinit'(V) \eqsp,
\end{multline*}
where $\rho_{\Cset}$ , $\Phi_{\Xinit,\Dset}$ and $\Psi_\Dset$ are defined in \eqref{eq:definition-rho}, \eqref{eq:definition-Phi} and \eqref{eq:definition-Psi}, respectively.
\end{cor}

\begin{pf}[Proof of Theorem \ref{thm:almost-sure-convergence}]
The conditions \eqref{eq:condition-Upsilon} and \eqref{eq:condition-Psi_D} imply that
\[
\limsup_{n \to \infty} \; \exp(-2 M n) \prod_{i=0}^n \Upsilon^2_{\Xset}(Y_i) \leq 1 \quad \text{and} \quad  \limsup_{n \to \infty} \; \exp(-2 Mn) \prod_{i=0}^n \Psi_\Dset^{-2}(Y_i) \leq 1 \eqsp.
\]
Condition (H\ref{assum:likelihood-not-zero}) and $\Xinit Q \1_\Dset > 0$
implies that $\Phi_{\Xinit,\Dset}(y_0,y_1) > 0$ for any $(y_0,y_1) \in \Yset^2$.
We then choose $\eta$ small enough so that
\[ \lim_{n \to \infty} \eta^{(\gamma-\beta) n/2} \exp(4Mn) (\epsilon_\Dset^-)^{-2(n-1)} = 0 \eqsp.
\]
The
proof follows from Corollary \ref{coro:main-statement}.
\end{pf}

\begin{pf}[Proof of Theorem \ref{thm:convergence-moment}]
Note that for any $\alpha \in (0,1)$ and any integer $n$,
\[
\tPE[ \tvnorm{\filt[\Xinit]{n}[\chunk{Y}{0}{n}] - \filt[\Xinit']{n}[\chunk{Y}{0}{n}]}] \leq \alpha^n +
\tPP[\tvnorm{\filt[\Xinit]{n}[\chunk{Y}{0}{n}] - \filt[\Xinit']{n}[\chunk{Y}{0}{n}]} \geq \alpha^n] \eqsp.
\]
Consider now the second term in the RHS of the previous equation.  Denote $\Omega_n$ the event
\begin{multline*}
  \Omega_n \eqdef \Bigl\{ \log \Phi_{\Xinit,\Dset}(Y_0,Y_1) > -M_0 n \eqsp,  \log \Phi_{\Xinit',\Dset} (Y_0,Y_1) >  -M_0 n \eqsp, \\
  \sum_{i=0}^n \log \Upsilon_\Xset(Y_i) < M_1 n, \sum_{i=2}^n \log
  \Psi_\Dset(Y_i) > - M_2 n, \sum_{i=1}^n \1_{\Kset}(Y_i) > n (1+\gamma)/2
  \Bigr\} \eqsp.
\end{multline*}
Clearly, $\tPP(\Omega^c_n) \leq \sum_{i=1}^3 r_i(n) + r_0(\Xinit,n) +
r_0(\Xinit',n)$ where $\{r_i(n)\}_{n \geq 0}$ and $\{ r_0(\Xinit,n) \}_{n \geq
  0}$ are defined in Eqs.~\eqref{eq:definition-r0}-\eqref{eq:definition-r3}. On
the event $\Omega_n$,
\[
\Phi^{-1}_{\Xinit,\Dset}(Y_0,Y_1)\; \Phi^{-1}_{\Xinit',\Dset}(Y_0,Y_1) \;
\prod_{i=2}^n \Psi^{-2}_\Dset (Y_i) \; \prod_{i=0}^n \Upsilon_{\Xset}(Y_i) \leq
\rme^{2 n \sum_{i=0}^2 M_i }\eqsp.
\]
One may choose $\eta > 0$ small enough and $\varrho \in (0,1)$ so that, for any $n$,
\[
\eta^{(\gamma-\beta)n/2}    \rme^{2 n \sum_{i=0}^2 M_i} (\epsilon_\Dset^-)^{-2(n-1)} \leq \varrho^n \eqsp.
\]
The proof then follows from Corollary \ref{coro:main-statement}.
\end{pf}

\section{Proof of Propositions \ref{prop:bound-condition-initiale}, \ref{prop:drift-condition-Polynomial}, and \ref{prop:drift-condition-MLP}}
\label{sec:ProofofpropositionsBoundConditionInitiale}

\begin{pf}[Proof of Proposition \ref{prop:bound-condition-initiale}]
  By the Jensen inequality with the function $u \mapsto [\log(u)]_-^p$, we
  obtain that for any $p\geq 1$,
\begin{multline}
\label{eq:upperboundlogPhi}
\left[ \log \Phi_{\Xinit,\Dset}(Y_0,Y_1) -  \log \left(\Xinit Q \1_\Dset\right) \right]_-^p \\ \leq
2^{p-1} \, (\Xinit Q \1_\Dset)^{-1} \iint \Xinit(dx_0) Q(x_0,dx_1) \1_\Dset(x_1) \; \sum_{i=0}^1 [ \log g(x_i,Y_i)]_-^p \eqsp,
\end{multline}
which implies by the Fubini theorem,
\begin{multline*}
  \tPE \left\{ \left[ \log \Phi_{\Xinit,\Dset}(Y_0,Y_1) - \log (\Xinit Q \1_\Dset) \right]^p_- \right\} \\
  \leq 2^{p-1} (\Xinit Q \1_\Dset)^{-1} \iint \Xinit(dx_0) Q(x_0,dx_1)
  \sum_{i=0}^1 \tPE [ \log g(x_i,Y_i)]^p_- \eqsp.
\end{multline*}
Since $\sup_{\Xset} V^{-1} \tPE[ \log g(\cdot,Y_i)]^p_- < \infty$, and $\sup_\Xset V^{-1} Q V < \infty$,
\begin{multline}
\label{eq:borne-Lp}
\iint \Xinit(dx_0) Q(x_0,dx_1) \sum_{i=0}^1 \tPE [ \log g(x_i,Y_i) ]^p_- \\
\leq \Xinit(V) \left\{ \sup_{i=0,1} \sup_{\Xset} V^{-1} \tPE[ \log g(\cdot,Y_i)]^p_- \right\} (1 + \sup_\Xset V^{-1} Q V) \eqsp.
\end{multline}
Similarly, for $\lambda > 0$, using the Jensen inequality with $u\mapsto \exp\left[ (\lambda/2) [\log u]_- \right]$ and the Fubini Theorem, we have
\begin{multline*}
\tPE\left[ \exp \left( (\lambda/2) [ \log \Phi_{\Xinit,\Dset}(Y_0,Y_1)  -  \log \left(\Xinit Q \1_\Dset\right)]_- \right) \right] \leq (\Xinit Q \1_\Dset)^{-1} \\
\times \iint \Xinit(dx_0) Q(x_0,dx_1) \tPE^{1/2}\left[ \exp(\lambda [\log g(x_0,Y_0)]_-) \right] \tPE^{1/2}\left[\exp(\lambda [\log g(x_1,Y_1)]_-)\right] \eqsp,
\end{multline*}
and the proof follows since $\sup_\Xset V^{-1/2} \Q V^{1/2} < \infty$.
\end{pf}

\begin{pf}[Proof of Proposition \ref{prop:drift-condition-Polynomial}]
  Let $\varphi$ be a non negative function on $\Yset$. Assume that
  $\sup_{\Xset} \tU^{-1} \tG(\cdot,\varphi^p) < \infty$. Proposition
  \ref{prop:drift-condition-LLN} shows that $\tStatDistrib \left[
    \tG(\cdot,\varphi^p) \right] < \infty$. Without loss of generality, we
  assume that $\tStatDistrib\left[ \tG(\cdot,\varphi) \right]= 0$.  For any $p
  \geq 1$,
\begin{multline}
  \PE^{\tJointKernel}_{\tXinit \otimes \tG} \left| \sum_{i=0}^n \varphi(Y_i)  \right|^p \\
  \leq 2^{p-1} \left( \PE^{\tJointKernel}_{\tXinit\otimes \tG } \left|
      \sum_{i=0}^n \{ \varphi(Y_i) - \tG(X_i,\varphi) \} \right|^p +
    \PE^{\tJointKernel}_{\tXinit \otimes \tG}\left| \sum_{i=0}^n
      \tG(X_i,\varphi) \right|^p \right) \eqsp.
\label{BlockY-BlockX}
\end{multline}
Since conditionally to $\chunk{X}{0}{n}$ the random variables $\chunk{Y}{0}{n}$ are
independent, we may apply the Marcinkiewicz-Zygmund inequality~\cite[Inequality
2.6.18 p.  82]{hall:heyde:1980}, showing that there exists a constant $c(p)$
depending only on $p$ such that
\begin{multline*}
  \PE^{\tJointKernel}_{\tXinit \otimes \tG} \left| \sum_{i=0}^n \{ \varphi(Y_i)
    - \tG(X_i,\varphi) \} \right|^p \leq c(p) \PE^{\tJointKernel}_{\tXinit
    \otimes \tG} \left( \sum_{i=0}^n |\varphi(Y_i) - \tG(X_i,\varphi)|^2
  \right)^{p/2} \eqsp.
\end{multline*}
If $1 \leq p \leq 2$,
\begin{multline*}
  \PE^{\tJointKernel}_{\tXinit \otimes \tG} \left| \sum_{i=0}^n \{ \varphi(Y_i)
    - \tG(X_i,\varphi) \} \right|^p \leq c(p) \sum_{i=0}^n
  \PE^{\tJointKernel}_{\tXinit \otimes \tG} \left[ |\varphi(Y_i) -
    \tG(X_i,\varphi)|^p
  \right] \\
  \leq 2^p c(p) \sum_{i=0}^n \PE^{\tJointKernel}_{\tXinit \otimes \tG} \left[
    \tG(X_i,|\varphi|^p) \right] \eqsp.
\end{multline*}
If $p \geq 2$, the Minkowski inequality yields
\begin{multline*}
  \PE^{\tJointKernel}_{\tXinit \otimes \tG} \left| \sum_{i=0}^n \{ \varphi(Y_i) -
    \tG(X_i,\varphi) \} \right|^p \leq c(p) \left( \sum_{i=0}^n
    \PE^{\tJointKernel}_{\tXinit \otimes \tG} \left[ |\varphi(Y_i) - \tG(X_i,\varphi)|^p
    \right]^{2/p} \right)^{p/2} \\
  \leq 2^p c(p) n^{p/2-1} \sum_{i=0}^n \PE^{\tJointKernel}_{\tXinit \otimes \tG} \left[
    \tG(X_i,|\varphi|^p) \right] \eqsp.
\end{multline*}
The $f$-norm ergodic theorem \cite[Theorem 14.0.1]{meyn:tweedie:1993} implies
that there exists a constant $C < \infty$, such that for any initial
probability measure $\tXinit$,
\[
\sum_{i=0}^n \PE^{\tJointKernel}_{\tXinit \otimes \tG} \left[
  \tG(X_i,|\varphi|^p) \right] \leq (n+1) \tStatDistrib \left(
  \tG\left(\cdot,\varphi^p \right) \right) + C \tXinit(\tV) \eqsp.
\]
Combining these discussions imply that there exists a finite constant $C_1$ such that
\[
\PE^{\tJointKernel}_{\tXinit \otimes \tG} \left| \sum_{i=0}^n \{ \varphi(Y_i) -
  \tG(X_i,\varphi) \} \right|^p \leq C_1 \; n^{p/2 \vee 1} \; \tXinit(\tV)
\eqsp.
\]
We now consider the second term in~(\ref{BlockY-BlockX}). Following the same
lines as in the proof of \cite[Proposition 12]{fort:moulines:2003} and applying
the Burkholder's inequality for martingales~\cite[Theorem
2.10]{hall:heyde:1980}, there exists a constant $C_2< \infty$ such that
\[
\PE^{\tJointKernel}_{\tXinit \otimes \tG}\left| \sum_{i=0}^n \tG(X_i,\varphi) \right|^p
\leq C_2 \; n^{p/2 \vee 1} \; \tXinit(\tW) \eqsp.
\]
The result follows.
\end{pf}

\begin{pf}[Proof of Proposition \ref{prop:drift-condition-MLP}]
  The first statement follows from standard results on phi-irreducible Markov
  chains satisfying the Foster-Lyapunov drift condition
  \cite{meyn:tweedie:1993}.  By Lemma \ref{lem:exponential-bound}, for any $x
  \in \Xset$ and $F \in \mathcal{G}_{\tW}$,
\begin{equation}
\label{eq:ExponentialInequality}
\PE_x^{\tQ} \left[ \exp \left( \sum_{k=0}^n F(X_k) \right) \right] \leq \tV(x) \rme^{(n+1)(\tb+\sup_\Xset (F-\tW))} \eqsp.
\end{equation}
Since under the probability $\PP_{\tXinit \otimes \tG}^{\tJointKernel}$ the
random variables $\chunk{Y}{0}{n}$ are conditionally independent given
$\chunk{X}{0}{n}$, and the conditional distribution of $Y_i$ given
$\chunk{X}{0}{n}$ is $\tG(X_i,\cdot)$,
\begin{multline*}
  \PE_{\tXinit \otimes \tG}^{\tJointKernel} \left[ \prod_{k=0}^n \rme^{\tlambda \varphi(Y_k)} \right] = \PE^{\tQ}_{\tXinit} \left[ \PE^{\tJointKernel} \left\{ \left. \prod_{k=0}^n \rme^{\tlambda \varphi(Y_k)} \right| \chunk{X}{0}{n} \right\}\right] \\
  = \PE_{\tXinit}^{\tQ} \left[ \prod_{k=0}^n \tG\left(X_k,\rme^{\tlambda \varphi}
    \right) \right] \leq \PE_{\tXinit}^{\tQ} \left[ \exp\left( \sum_{k=0}^n
      \left| \log \tG\left(X_k,\rme^{\tlambda \varphi} \right) \right| \right)
  \right] \eqsp.
\end{multline*}
By the Jensen inequality, $F \eqdef \log \tG \left( \cdot, \rme^{\tlambda \varphi}
\right)$ is non negative and belongs to $\mathcal{G}_{\tW}$; we may thus apply
\eqref{eq:ExponentialInequality} which yields
\[
\PE_{\tXinit  \otimes \tG}^{\tJointKernel} \left[ \prod_{k=0}^n \rme^{\tlambda \varphi(Y_k)} \right] \leq \tXinit \left( \tV \right) \rme^{(n+1)(\tb+\sup_\Xset (F-\tW))} \eqsp.
\]
The proof then follows by applying the Markov inequality.
\end{pf}

\appendix
\section{Technical Results}
\label{sec:Technical-Results}
We have collected in this section the proof of some of the technical results.
\begin{lem}
\label{lem:denombrement}
For any integer $n \geq 1$, and sequence $\xb \eqdef \{ x_i \}_{i \geq 0} \in \{0,1\}^{\nset}$, denote by $M_n(\xb) \eqdef \sum_{i=0}^{n-1} \1 \{ x_i=1 \}$ and
$N_n(\xb) \eqdef \sum_{i=0}^{n-1} \1 \{x_i=1, x_{i+1}=1 \}$. Then,
\[
M_n(\xb) \leq \frac{n+1}{2} + \frac{N_n(\xb)}{2} \eqsp.
\]
\end{lem}
\begin{pf}
Denote by $\tau$ the shift operator on sequences defined, for any sequence $\xb \eqdef \{ x_i \}_{i \geq 1}$, by $[\tau \xb]_k = x_{k+1}$.
Let $\mathbf{x}= \{ x_i \}_{i \geq 0}$ be a sequence such that $x_j= 0$ for $j \geq n$.
By construction, $N_n(\xb)= M_n( \xb \, \mathrm{AND} \, \tau \xb)$. The proof then follows from the obvious identity:
\begin{multline*}
n \geq M_{n}( \xb \, \mathrm{OR}  \, \tau \xb) = M_n( \xb) + M_n (\tau \xb) - M_n( \xb \, \mathrm{AND} \, \tau \xb) \\
\geq 2 M_n(\xb) - 1 - N_n(\xb) \eqsp,
\end{multline*}
where $\mathrm{AND}$ and $\mathrm{OR}$ is the componentwise incluse "AND" and "OR".
\end{pf}

\begin{lem}
\label{lem:exponential-bound}
Assume that there exist a function $V: \Xset \to [1,\infty)$, a function $W: \Xset \to (0,\infty)$  and a constant $b < \infty$
such that
\begin{equation}
\label{eq:drift-condition-MLP-1}
\log( V^{-1} Q V) \leq - W + b \eqsp.
\end{equation}
Let $n$ be an integer and $F_k$, $k=0,\dots,n-1$, be functions belonging to $\mathcal{G}_W$, where
$\mathcal{G}_{W}$ is defined in \eqref{eq:definition-mathcalG}. Hence, for any $x \in \Xset$,
\begin{equation}
\label{eq:ExponentialInequality-0}
\PE_x^{Q} \left[ \exp \left( \sum_{k=0}^{n-1} |F_k(X_k)| \right) \right] \leq V(x) \rme^{bn + \sum_{k=0}^{n-1} \sup_\Xset (|F_k| - W)} \eqsp.
\end{equation}
\end{lem}
\begin{pf}
The proof is adapted from \cite[Theorem 2.1]{kontoyiannis:meyn:2005}. Set for any integer $n$,
\begin{equation}
\label{eq:SuperMartingale}
M_n \eqdef  V(X_n) \exp \left( \sum_{k=0}^{n-1} \left\{  W(X_k) - b  \right\} \right) \eqsp.
\end{equation}
The multiplicative drift condition \eqref{eq:drift-condition-MLP-1} implies that $\{ M_n \}$ is a supermartingale. Hence, for any $n \in \nset$ and $x \in \Xset$,
\[
\PE_x^{Q} \left[ V(X_n) \exp \left(  - b n + \sum_{k=0}^{n-1}  W(X_k) \right) \right] \leq V(x) \eqsp.
\]
The proof follows.
\end{pf}


\begin{thebibliography}{10}

\bibitem{andrieu:doucet:2002}
C.~Andrieu and A.~Doucet.
\newblock Particle filtering for partially observed {G}aussian state space
  models.
\newblock {\em J. Roy. Statist. Soc. Ser. B}, 64(4):827--836, 2002.

\bibitem{atar:zeitouni:1997}
R.~Atar and O.~Zeitouni.
\newblock Exponential stability for nonlinear filtering.
\newblock {\em Ann. Inst. H. Poincar\'e Probab. Statist.}, 33(6):697--725,
  1997.

\bibitem{budhiraja:ocone:1997}
A.~Budhiraja and D.~Ocone.
\newblock Exponential stability of discrete-time filters for bounded
  observation noise.
\newblock {\em Systems Control Lett.}, 30:185--193, 1997.

\bibitem{budhiraja:ocone:1999}
A.~Budhiraja and D.~Ocone.
\newblock Exponential stability in discrete-time filtering for non-ergodic
  signals.
\newblock {\em Stochastic Process. Appl.}, 82(2):245--257, 1999.

\bibitem{cappe:moulines:ryden:2005}
O.~Capp\'{e}, E.~Moulines, and T.~Ryd\'{e}n.
\newblock {\em Inference in Hidden {M}arkov Models}.
\newblock Springer, 2005.

\bibitem{chigansky:lipster:2004}
P.~Chigansky and R.~Lipster.
\newblock Stability of nonlinear filters in nonmixing case.
\newblock {\em Ann. Appl. Probab.}, 14(4):2038--2056, 2004.

\bibitem{dedecker:doukhan:2003}
J.~Dedecker and P.~Doukhan.
\newblock A new covariance inequality and applications.
\newblock {\em Stochastic Process. Appl.}, 106(1):63--80, 2003.

\bibitem{delmoral:2004}
P.~{Del Moral}.
\newblock {\em {F}eynman-Kac {F}ormulae. {G}enealogical and Interacting
  Particle Systems with Applications}.
\newblock Springer, 2004.

\bibitem{delmoral:guionnet:1998}
P.~{Del Moral} and A.~Guionnet.
\newblock Large deviations for interacting particle systems: applications to
  non-linear filtering.
\newblock {\em Stoch. Proc. App.}, 78:69--95, 1998.

\bibitem{delmoral:guionnet:2001}
P.~{Del Moral} and A.~Guionnet.
\newblock On the stability of interacting processes with applications to
  filtering and genetic algorithms.
\newblock {\em Annales de l'Institut Henri Poincar\'e}, 37:155--194, 2001.

\bibitem{delmoral:ledoux:miclo:2003}
P.~Del~Moral, M.~Ledoux, and L.~Miclo.
\newblock On contraction properties of {M}arkov kernels.
\newblock {\em Probab. Theory Related Fields}, 126(3):395--420, 2003.

\bibitem{fort:moulines:2003}
G.~Fort and E.~Moulines.
\newblock Convergence of the {M}onte {C}arlo expectation maximization for
  curved exponential families.
\newblock {\em Ann. Statist.}, 31(4):1220--1259, 2003.

\bibitem{hall:heyde:1980}
P.~Hall and C.~C. Heyde.
\newblock {\em Martingale Limit Theory and its Application}.
\newblock Academic Press, New York, London, 1980.

\bibitem{hull:white:1987}
J.~Hull and A.~White.
\newblock The pricing of options on assets with stochastic volatilities.
\newblock {\em J. Finance}, 42:281--300, 1987.

\bibitem{jacquier:polson:rossi:1994}
E.~Jacquier, N.~G. Polson, and P.~E. Rossi.
\newblock {B}ayesian analysis of stochastic volatility models (with
  discussion).
\newblock {\em J. Bus. Econom. Statist.}, 12:371--417, 1994.

\bibitem{kleptsyna:veretennikov:2007}
M.L. Kleptsyna and A.Y. Veretennikov.
\newblock On discrete time ergodic filters with wrong initial conditions.
\newblock Technical report, Université du Maine and University of Leeds, 2007.
\newblock Available at
  http://www.univ-lemans.fr/sciences/statist/download/Kleptsyna/filt19fg.pdf.

\bibitem{kontoyiannis:meyn:2005}
I.~Kontoyiannis and S.~P. Meyn.
\newblock Large deviations asymptotics and the spectral theory of
  multiplicatively regular {M}arkov processes.
\newblock {\em Electron. J. Probab.}, 10:no. 3, 61--123 (electronic), 2005.

\bibitem{legland:oudjane:2003}
Fran{\c{c}}ois LeGland and Nadia Oudjane.
\newblock A robustification approach to stability and to uniform particle
  approximation of nonlinear filters: the example of pseudo-mixing signals.
\newblock {\em Stochastic Process. Appl.}, 106(2):279--316, 2003.

\bibitem{manrique:shephard:1998}
A.~Manrique and N.~Shephard.
\newblock Likelihood inference for limited dependent processes.
\newblock {\em Econometrics Journal}, 1:174--202, 1998.

\bibitem{meyn:tweedie:1993}
S.~P. Meyn and R.~L. Tweedie.
\newblock {\em {M}arkov Chains and Stochastic Stability}.
\newblock Springer, London, 1993.

\bibitem{ocone:pardoux:1996}
D.~Ocone and E.~Pardoux.
\newblock Asymptotic stability of the optimal filter with respect to its
  initial condition.
\newblock {\em SIAM J. Control}, 34:226--243, 1996.

\bibitem{oudjane:rubenthaler:2005}
Nadia Oudjane and Sylvain Rubenthaler.
\newblock Stability and uniform particle approximation of nonlinear filters in
  case of non ergodic signals.
\newblock {\em Stoch. Anal. Appl.}, 23(3):421--448, 2005.

\end{thebibliography}
\end{document}